\input amstex
\let\myfrac=\frac
\input eplain.tex
\let\frac=\myfrac
\input epsf




\loadeufm \loadmsam \loadmsbm
\message{symbol names}\UseAMSsymbols\message{,}

\font\myfontdefault=cmr10

\font\mytdmchapfont=cmb10 at 14pt
\font\mytdmheadfont=cmb10 at 10pt
\font\mytdmsubheadfont=cmr10

\magnification 1200
\newif\ifinappendices
\newif\ifundefinedreferences
\newif\ifchangedreferences
\newif\ifloadreferences
\newif\ifmakebiblio
\newif\ifmaketdm

\undefinedreferencesfalse
\changedreferencesfalse


\loadreferencestrue
\makebibliofalse
\maketdmfalse

\def\headpenalty{-400}     
\def\proclaimpenalty{-200} 

%
%

\def\alphanum#1{\ifcase #1 _\or A\or B\or C\or D\or E\or F\or G\or H\or I\or J\or K\or L\or M\or N\or O\or P\or Q\or R\or S\or T\or U\or V\or W\or X\or Y\or Z\fi}
\def\gobbleeight#1#2#3#4#5#6#7#8{}

\newwrite\references
\newwrite\tdm
\newwrite\biblio

\newcount\chapno
\newcount\headno
\newcount\subheadno
\newcount\procno
\newcount\figno
\newcount\citationno

\def\setcatcodes{%
\catcode`\!=0 \catcode`\\=11}%

\ifloadreferences
    {\catcode`\@=11 \catcode`\_=11%
    \global\def\_@citation@Cheeger{1}
\global\def\_@citation@Corlette{2}
\global\def\_@citation@GromA{3}
\global\def\_@citation@Peterson{4}
\global\def\_@citation@SmiA{5}
\global\def\_@citation@SmiE{6}
\global\def\_@proc@DefnHigherDerivativesSecondFF{1.1}
\global\def\_@proc@TheoremTC{1.2}
\global\def\_@proc@LemmaBDONGBDONDPHI{2.3}
\global\def\_@proc@LemmaBDONDKGBDONDKPHI{2.4}
\global\def\_@proc@CorUniquenessOfCheegerGromov{2.5}
\global\def\_@proc@LemmaUniquenessOfGraphFunction{3.1}
\global\def\_@proc@LemmaDisasterConditions{3.3}
\global\def\_@proc@LemmaControlOfExponential{3.4}
\global\def\_@proc@LemmaTransformationsOfBoundsOnII{3.5}
\global\def\_@proc@LemmaBoundIIEuc{3.6}
\global\def\_@proc@LemmaBoundOnDDf{3.7}
\global\def\_@proc@LemmaGraphFunctionsOverGivenRadius{3.8}
\global\def\_@proc@LemmaControlZeroOrder{3.9}
\global\def\_@proc@LemmaProjOntoTGraph{3.10}
\global\def\_@proc@CorSizeOfProjOntoTGraph{3.11}
\global\def\_@proc@LemmaControlAllOrders{3.12}
\global\def\_@proc@LemmaControlInTermsOfII{3.13}
\global\def\_@proc@LemmaBDONTRANSITION{4.2}
\global\def\_@proc@LemmaLA{4.3}
\global\def\_@proc@LemmaLB{4.4}
\global\def\_@head@HolderSpaces{A}
\global\def\_@subhead@ManifoldConditions{A.1}
\global\def\_@proc@LemmaBasicLipschitzComposition{A.1}
\global\def\_@proc@LemmaAdvancedLipschitzComposition{A.2}
\global\def\_@proc@LemmaBasicCtyAddMult{A.4}
\global\def\_@proc@LemmaCtyAddMult{A.5}
\global\def\_@proc@LemmaBasicCtyComposition{A.7}
\global\def\_@proc@LemmaCtsComp{A.9}
\global\def\_@proc@DefnOptimalAtlas{B.1}
\global\def\_@proc@ThmConvergenceOfRiemannianGeometry{B.3}
\global\def\_@proc@PropDiffPlusRiemStruct{B.7}
\global\def\_@proc@PropConvergenceMappings{B.9}
    }%
\else
    \openout\references=references.tex
\fi

\newcount\newchapflag 
\newcount\showpagenumflag 

\global\chapno = -1 
\global\citationno=0
\global\headno = 0
\global\subheadno = 0
\global\procno = 0
\global\figno = 0

\def\resetcounters{%
\global\headno = 0%
\global\subheadno = 0%
\global\procno = 0%
\global\figno = 0%
}

\global\newchapflag=0 
\global\showpagenumflag=0 

\def\headticket{\ifinappendices\alphanum\headno\else\the\headno\fi}%

\def\chinfo{\ifinappendices\alphanum\chapno\else\the\chapno\fi}%
\def\headinfo{\headticket}
\def\subheadinfo{\headticket.\the\subheadno}
\def\procinfo{\headticket.\the\procno}
\def\figinfo{\headticket.\the\figno}
\def\citationinfo{\the\citationno}%
\def\nextheadno{\global\advance\headno by 1 \global\subheadno = 0 \global\procno = 0}
\def\nextsubheadno{\global\advance\subheadno by 1}
\def\nextprocno{\global\advance\procno by 1 \procinfo}
\def\nextfigno{\global\advance\figno by 1 \figinfo}

{\global\let\noe=\noexpand%
%
%
\catcode`\@=11%
\catcode`\_=11%
\setcatcodes%
!global!def!_@@internal@@makeref#1{%
!global!expandafter!def!csname #1ref!endcsname##1{%
!csname _@#1@##1!endcsname%
!expandafter!ifx!csname _@#1@##1!endcsname!relax%
    !write16{#1 ##1 not defined - run saving references}%
    !undefinedreferencestrue%
!fi}}%
!global!def!_@@internal@@makelabel#1{%
!global!expandafter!def!csname #1label!endcsname##1{%
!edef!temptoken{!csname #1info!endcsname}%
!ifloadreferences%
    !expandafter!ifx!csname _@#1@##1!endcsname!relax%
        !write16{#1 ##1 not hitherto defined - rerun saving references}%
        !changedreferencestrue%
    !else%
        !expandafter!ifx!csname _@#1@##1!endcsname!temptoken%
        !else
            !write16{#1 ##1 reference has changed - rerun saving references}%
            !changedreferencestrue%
        !fi%
    !fi%
!else%
    !expandafter!edef!csname _@#1@##1!endcsname{!temptoken}%
    !edef!textoutput{!write!references{\global\def\_@#1@##1{!temptoken}}}%
    !textoutput%
!fi}}%
!global!def!makecounter#1{!_@@internal@@makelabel{#1}!_@@internal@@makeref{#1}}%
!unsetcatcodes%
}
\makecounter{ch}%
\makecounter{head}%
\makecounter{subhead}%
\makecounter{proc}%
\makecounter{fig}%
\makecounter{citation}%
\def\newref#1#2{%
\def\temptext{#2}%
\edef\bibliotextoutput{\expandafter\gobbleeight\meaning\temptext}%
\global\advance\citationno by 1\citationlabel{#1}%
\ifmakebiblio%
    \edef\fileoutput{\write\biblio{\noindent\hbox to 0pt{\hss$[\the\citationno]$}\hskip 0.2em\bibliotextoutput\medskip}}%
    \fileoutput%
\fi}%
\def\cite#1{%
$[\citationref{#1}]$%
\ifmakebiblio%
    \edef\fileoutput{\write\biblio{#1}}%
    \fileoutput%
\fi%
}%
%
%
%

\let\mypar=\par


\def\raggedleft{\leftskip=0pt plus 1fil \parfillskip=0pt}


\font\lettrinefont=cmr10 at 28pt
\def\lettrine #1[#2][#3]#4%
{\hangafter -#1 \hangindent #2
\noindent\hskip -#2 \vtop to 0pt{
\kern #3 \hbox to #2 {\lettrinefont #4\hss}\vss}}

\font\mylettrinefont=cmr10 at 28pt
\def\mylettrine #1[#2][#3][#4]#5%
{\hangafter -#1 \hangindent #2
\noindent\hskip -#2 \vtop to 0pt{
\kern #3 \hbox to #2 {\mylettrinefont #5\hss}\vss}}


\edef\Pagetitle={Blank}

\headline={\hfil\Pagetitle\hfil}

\footline={\hfil\myfontdefault\folio\hfil}

\def\nextoddpage
{
\newpage%
\ifodd\pageno%
\else%
    \global\showpagenumflag = 0%
    \null%
    \vfil%
    \eject%
    \global\showpagenumflag = 1%
\fi%
}


\def\newchap#1#2%
{%
%
%
\global\advance\chapno by 1%
\resetcounters%
%
%
\newpage%
\ifodd\pageno%
\else%
    \global\showpagenumflag = 0%
    \null%
    \vfil%
    \eject%
    \global\showpagenumflag = 1%
\fi%
\global\newchapflag = 1%
\global\showpagenumflag = 1%
%
%
{\font\chapfontA=cmsl10 at 30pt%
\font\chapfontB=cmsl10 at 25pt%
\null\vskip 5cm%
{\chapfontA\raggedleft\hfil%
{%
\ifnum\chapno=0
    \phantom{%
    \ifinappendices%
        Annexe \alphanum\chapno%
    \else%
        \the\chapno%
    \fi}%
\else%
    \ifinappendices%
        Annexe \alphanum\chapno%
    \else%
        \the\chapno%
    \fi%
\fi%
}%
\par}%
\vskip 2cm%
{\chapfontB\raggedleft%
\lineskiplimit=0pt%
\lineskip=0.8ex%
\hfil #1\par}%
\vskip 2cm%
}%
\edef\Pagetitle{#2}%
%
%
\ifmaketdm%
    \def\temp{#2}%
    \def\tempbis{\nobreak}%
    \edef\chaptitle{\expandafter\gobbleeight\meaning\temp}%
    \edef\mynobreak{\expandafter\gobbleeight\meaning\tempbis}%
    \edef\textoutput{\write\tdm{\bigskip{\noexpand\mytdmchapfont\noindent\chinfo\ - \chaptitle\hfill\noexpand\folio}\par\mynobreak}}%
\fi%
\textoutput%
}


\def\newhead#1%
{%
\ifhmode%
    \mypar%
\fi%
\ifnum\headno=0%
    \ifinappendices%
        \nobreak\vskip -\lastskip%
        \nobreak\vskip .5cm%
        \line{\hfil$\diamond$\hfil}%
        \penalty\headpenalty\vskip .5cm%
    \fi
\else%
    \nobreak\vskip -\lastskip%
    \nobreak\vskip .5cm%
    \line{\hfil$\diamond$\hfil}%
    \penalty\headpenalty\vskip .5cm%
\fi%
\nextheadno%
\ifmaketdm%
    \def\temp{#1}%
    \edef\sectiontitle{\expandafter\gobbleeight\meaning\temp}%
    \edef\textoutput{\write\tdm{\noindent{\noexpand\mytdmheadfont\quad\headinfo\ - \sectiontitle\hfill\noexpand\folio}\par}}%
    \textoutput%
\fi%
\font\headfontA=cmbx10 at 14pt%
{\headfontA\noindent\headinfo\ -\ #1.\hfil}%
\nobreak\vskip .5cm%
}%


\def\newsubhead#1%
{%
\ifhmode%
    \mypar%
\fi%
\ifnum\subheadno=0%
\else%
    \penalty\headpenalty\vskip .4cm%
\fi%
\nextsubheadno%
\ifmaketdm%
    \def\temp{#1}%
    \edef\subsectiontitle{\expandafter\gobbleeight\meaning\temp}%
    \edef\textoutput{\write\tdm{\noindent{\noexpand\mytdmsubheadfont\quad\quad\subheadinfo\ - \subsectiontitle\hfill\noexpand\folio}\par}}%
    \textoutput%
\fi%
\font\subheadfontA=cmsl10 at 12pt
{\subheadfontA\noindent\subheadinfo\ #1.\hfil}%
\nobreak\vskip .25cm%
}%

%
%


\font\mathromanten=cmr10
\font\mathromanseven=cmr7
\font\mathromanfive=cmr5
\newfam\mathromanfam
\textfont\mathromanfam=\mathromanten
\scriptfont\mathromanfam=\mathromanseven
\scriptscriptfont\mathromanfam=\mathromanfive
\def\mathroman{\fam\mathromanfam}


\font\sansseriften=cmss10
\font\sansserifseven=cmss7
\font\sansseriffive=cmss5
\newfam\sansseriffam
\textfont\sansseriffam=\sansseriften
\scriptfont\sansseriffam=\sansserifseven
\scriptscriptfont\sansseriffam=\sansseriffive
\def\mathsf{\fam\sansseriffam}


\font\boldten=cmb10
\font\boldseven=cmb7
\font\boldfive=cmb5
\newfam\mathboldfam
\textfont\mathboldfam=\boldten
\scriptfont\mathboldfam=\boldseven
\scriptscriptfont\mathboldfam=\boldfive
\def\mathbf{\fam\mathboldfam}


\font\mycmmiten=cmmi10
\font\mycmmiseven=cmmi7
\font\mycmmifive=cmmi5
\newfam\mycmmifam
\textfont\mycmmifam=\mycmmiten
\scriptfont\mycmmifam=\mycmmiseven
\scriptscriptfont\mycmmifam=\mycmmifive

\def\hexa#1{\ifcase #1 0\or 1\or 2\or 3\or 4\or 5\or 6\or 7\or 8\or 9\or A\or B\or C\or D\or E\or F\fi}
\mathchardef\mathi="7\hexa\mycmmifam7B
\mathchardef\mathj="7\hexa\mycmmifam7C


\font\mymsbmten=msbm10 at 8pt
\font\mymsbmseven=msbm7 at 5.6pt
\font\mymsbmfive=msbm5 at 4pt
\newfam\mymsbmfam
\textfont\mymsbmfam=\mymsbmten
\scriptfont\mymsbmfam=\mymsbmseven
\scriptscriptfont\mymsbmfam=\mymsbmfive

\mathchardef\mybeth="7\hexa\mymsbmfam69
\mathchardef\mygimmel="7\hexa\mymsbmfam6A
\mathchardef\mydaleth="7\hexa\mymsbmfam6B


\def\placelabel[#1][#2]#3{{%
\setbox10=\hbox{\raise #2cm \hbox{\hskip #1cm #3}}%
\ht10=0pt%
\dp10=0pt%
\wd10=0pt%
\box10}}%


\newif\ifinproclaim%
\global\inproclaimfalse%
\def\proclaim#1{%
\medskip%
%
%
\bgroup%
\inproclaimtrue%
\setbox10=\vbox\bgroup\leftskip=0.8em\noindent{\bf #1}\sl%
}

\def\endproclaim{%
\egroup%
\setbox11=\vtop{\noindent\vrule height \ht10 depth \dp10 width 0.1em}%
\wd11=0pt%
\setbox12=\hbox{\copy11\kern 0.3em\copy11\kern 0.3em}%
\wd12=0pt%
\setbox13=\hbox{\noindent\box12\box10}%
\noindent\unhbox13%
\egroup%
\medskip\ignorespaces%
}

\def\proclaim#1{%
\medskip%
\bgroup%
\inproclaimtrue%
\noindent{\bf #1}%
\nobreak\medskip%
\sl%
}

\def\endproclaim{%
\mypar\egroup\penalty\proclaimpenalty\medskip\ignorespaces%
}

\def\noskipproclaim#1{%
\medskip%
\bgroup%
\inproclaimtrue%
\noindent{\bf #1}\nobreak\sl%
}

\def\endnoskipproclaim{%
\mypar\egroup\penalty\proclaimpenalty\medskip\ignorespaces%
}


\def\ninn{{n\in\Bbb{N}}}

\def\colon{\hphantom{0}:}
\def\proof{{\noindent\bf Proof:\ }}

\def\remark{{\noindent\sl Remark:\ }}
\def\suite#1#2{({#1}_{#2})_{#2\in\Bbb{N}}}

\def\msup{\mathop{{\mathroman Sup}}}
\def\minf{\mathop{{\mathroman Inf}}}
\def\msf#1{{\mathsf #1}}

\def\qed{~$\square$}
\def\munion{\mathop{\cup}}
\def\minter{\mathop{\cap}}
\def\myitem#1{%
\noindent\hbox to .5cm{\hfill#1\hss}
}%

\catcode`\@=11
\def\Eqalign#1{\null\,\vcenter{\openup\jot\m@th\ialign{%
\strut\hfil$\displaystyle{##}$&$\displaystyle{{}##}$\hfil%
&&\quad\strut\hfil$\displaystyle{##}$&$\displaystyle{{}##}$%
\hfil\crcr #1\crcr}}\,}
\catcode`\@=12

\def\makeop#1{%
\global\expandafter\def\csname op#1\endcsname{{\mathroman #1}}}%

\def\makeopsmall#1{%
\global\expandafter\def\csname op#1\endcsname{{\mathroman{\lowercase{#1}}}}}%

\makeopsmall{ArcTan}%
\makeopsmall{ArcCos}%
\makeop{Sup}%
\makeop{Arg}%
\makeop{Det}%
\makeop{Log}%
\makeop{Re}%
\makeop{Im}%
\makeop{Dim}%
\makeopsmall{Tan}%
\makeop{Cov}%
\makeop{Ker}%
\makeopsmall{Cos}%
\makeopsmall{Sin}%
\makeop{Exp}%
\makeopsmall{Tanh}%
\makeop{Tr}%
\makeop{End}%
\makeop{Long}%
\makeop{Ch}%
\makeop{Exp}%
\makeop{Int}%
\makeop{Ext}%
\makeop{Aire}%
\makeop{Length}%
\makeop{Im}%
\makeop{Conf}%
\makeop{Exp}%
\makeop{Mod}%
\makeop{Log}%
\makeop{Ext}%
\makeop{Int}%
\makeop{Dist}%
\makeop{Aut}%
\makeop{Id}%
\makeop{SO}%
\makeop{Homeo}%
\makeop{Vol}%
\makeop{Ric}%
\makeop{Hess}%
\makeop{Euc}%
\makeop{Isom}%
\makeop{Max}%
\makeop{Long}%
\makeop{Fixe}%
\makeop{Wind}%
\makeop{Mush}%
\makeop{Supp}%
\makeop{Ad}%
\makeop{loc}%
\makeop{Len}%
\makeop{Area}%
\makeop{SL}%
\makeop{GL}%
\makeop{dVol}%
\makeop{Min}%
\makeop{Symm}%
\makeop{Diam}%
\makeop{Lip}%
\makeop{Dil}%
\makeop{O}%
\makeop{Rest}%
\let\emph=\bf

\hyphenation{quasi-con-formal}

%
%

\ifmakebiblio%
    \openout\biblio=biblio.tex %
    {%
        \edef\fileoutput{\write\biblio{\bgroup\leftskip=2em}}%
        \fileoutput
    }%
\fi%

\newref{Cheeger}{Cheeger J., Finiteness theorems for Riemannian manifolds, {\sl Amer. J. Math.} {\bf 92} (1970), 61--74}
\newref{Corlette}{Corlette K., Immersions with bounded curvature, {\sl Geom. Dedicata} {\bf 33} (1990), no. 2, 153--161}
\newref{GromA}{Gromov M., {\sl Metric Structures for Riemannian and Non-Riemannian Spaces}, Progress in Mathematics, {\bf 152}, Birkh\"auser, Boston, (1998)}
\newref{Peterson}{Peterson P., {\sl Riemannian Geometry}, Graduate Texts in Mathematics, {\bf 171}, Springer Verlag, New York, (1998)}
\newref{SmiA}{Smith G., Special Legendrian structures and Weingarten problems, Preprint, Orsay (2005)}%
\newref{SmiE}{Smith G., Th\`ese de doctorat, Paris (2004)}%

\ifmakebiblio%
    {\edef\fileoutput{\write\biblio{\egroup}}%
    \fileoutput}%
\fi%

%
%
%
\document
\inappendicesfalse
\myfontdefault
\global\chapno=1
\global\showpagenumflag=1
\def\Pagetitle{}
\null
\vfill
\def\centre{\rightskip=0pt plus 1fil \leftskip=0pt plus 1fil \spaceskip=.3333em \xspaceskip=.5em \parfillskip=0em \parindent=0em}%
\def\textmonth#1{\ifcase#1\or January\or Febuary\or March\or April\or May\or June\or July\or August\or September\or October\or November\or December\fi}
\font\abstracttitlefont=cmr10 at 14pt
{\abstracttitlefont\centre An Arzela-Ascoli Theorem for Immersed Submanifolds\par}
\bigskip
{\centre Graham Smith\par}
\bigskip
{\centre \the\day\ \textmonth\month\ \the\year\par}
\bigskip
{\centre Equipe de topologie et dynamique,\par
Laboratoire de math\'ematiques,\par
B\^atiment 425,\par
UFR des sciences d'Orsay,\par
91405 Orsay CEDEX, FRANCE\par
UMR 8628 du CNRS\par}
\bigskip
\noindent{\emph Abstract:\ }The classical Arzela-Ascoli theorem is a compactness result for families of functions depending on bounds on the derivatives of the functions, and is of invaluable use in many fields of mathematics. In this paper, inspired by a result of Corlette, we prove an analogous compactness result for families of immersed submanifolds which depends only on bounds on the derivatives of the second fundamental forms of these submanifolds. We then show how 
the result of Corlette may be obtained as an immediate corollary.
\par
\vfill
\nextoddpage
\def\Pagetitle{\sl An Arzela-Ascoli Theorem for Immersed Submanifolds}
\def\lipa{\text{Lip}^\alpha}

\def\lip#1{\text{Lip}^{#1}}
\def\slipa{\text{Lip}^\alpha(\Omega)}
\def\lka{\text{C}^{k,\alpha}}
\def\lkaloc{\lka_\oploc}
\def\slka{\text{C}^{k,\alpha}(\Omega)}
\def\slkaloc{\text{C}^{k,\alpha}_\oploc(\Omega)}

\def\nslipa#1{\|#1\|_{\lipa(\Omega)}}
\def\nsczero#1{\|#1\|_{C^0(\Omega)}}

\def\nslka#1{\|#1\|_{C^{k,\alpha}(\Omega)}}

\newhead{Introduction}
\noindent In \cite{Cheeger}, Cheeger proved his famous ``Finiteness Theorem'' which states that,
given positive real numbers $K,\epsilon,R\in\Bbb{R}^+$, there exist only finitely many
homeomorphism classes of complete manifolds of a given dimension with sectional curvature
bounded above in absolute value by $K$, injectivity radius bounded below by $\epsilon$ and
diameter bounded above by $R$. Gromov later showed how this result may be viewed in terms of
(more precisely, as a corollary of) a compactness result for the family of such manifolds (see,
for example \cite{GromA} and \cite{Peterson}). Viewed from this perspective, the conditions imposed by Cheeger become quite satisfying to our common sense. The curvature bound is probably the most fundamental. Indeed, curvature reflects the ``derivatives'' of the manifold (how fast it turns), and whenever derivatives are bounded, following the philosophy of the classical Arzela-Ascoli theorem, one is entitled to expect to find a compactness result. The other two conditions reflect more geometrical considerations. The lower bound on the injectivity radius excludes ``pinching'' (one can consider a sequence of cylinders of ever smaller radius: the only intrinsic data that informs us of degeneration is the injectivity radius which tends to zero), and the bound on the diameter excludes the possibility of adding components indefinitely without introducing very high curvatures (to see what happens without this condition, one can consider a sequence of surfaces of ever increasing genus).
\medskip
\noindent In \cite{Corlette}, Corlette proved an analogous finiteness theorem for immersed submanifolds of a given Riemannian manifold. He proves that, given positive real numbers 
$K,R\in\Bbb{R}^+$ and a given compact manifold $M$, there exists only finitely many $C^1$ 
isotopy classes of complete immersed submanifolds of $M$ with second fundamental form bounded
above in norm by $K$ and diameter bounded above by $R$. This result no longer requires the
condition on the injectivity radius of the immersed submanifold since such a lower bound is now a product of the bounds on the second fundamental form of the submanifold and the curvature of the ambient manifold (one may consider again the example of a cylinder in Euclidean space: its injectivity radius cannot become small without its second fundamental form becoming large). Following the philosophy of Cheeger's finiteness theorem, one would expect this result to also arise from a compactness result for immersed submanifolds.
\medskip
\noindent The compactness result thus obtained bears a perfect analogy to the classical Arzela-Ascoli theorem, being, in certain aspects, a generalisation of this result, and it is for this reason that we have chosen to name it thus. An
example of an application of this result may be found in \cite{SmiA}. The statement of this theorem requires the following definition:
\proclaim{Definition \nextprocno}
\noindent Let $(M,g)$ be a Riemannian manifold. Let $X=(Y,i)$ be an immersed submanifold in $M$. Let $\nabla^i$ be the Levi-Civita covariant derivative generated over $Y$ by
the immersion $i$ into $(M,g)$. Let $A(X)$ be the second fundamental form of $X$. For all $k\geqslant 2$, we define $A_k(X)$ using the following recurrence relation:
$$\matrix
A_2(X) \hfill&= A(X), \hfill\cr
A_k(X) \hfill&= \nabla^i A_{k-1}(X)\ \forall k\geqslant 3.\hfill\cr
\endmatrix$$
\noindent We now define $\Cal{A}_k(X)$ for all $k\geqslant 2$ by:
$$
\Cal{A}_k(X) = \sum_{i=2}^k \|A_i(X)\|.
$$
\endproclaim
\proclabel{DefnHigherDerivativesSecondFF}
\noindent The principal result of this paper is now given by the following theorem:
\proclaim{Theorem \nextprocno}
\noindent Let $k\geqslant 2, m\leqslant n\in\Bbb{N}$ be positive integers.
\medskip
\noindent Let $(M_n,p_n)_\ninn,(M_\infty,p_\infty)$ be complete pointed Riemannian manifolds of dimension $n$ and of class at least
$C^k$ such that $(M_n,p_n)_\ninn$ converges towards $(M_\infty,p_\infty)$ in the pointed $C^k$ Cheeger/Gromov topology.
\medskip
\noindent For all $n\in\Bbb{N}$, let $\Sigma_n=(S_n,q_n)$ be an $m$-dimensional pointed immersed submanifold of $M_n$ of type at least $C^k$ such that $i(q_n)=p_n$.
\medskip
\noindent Suppose that for all $R>0$, there exists $B$ such that, for all $n$:
$$
\left|\Cal{A}_{\Sigma_n}^k(q)\right|\leqslant B\qquad \forall q\in B_R(q_n).
$$
\noindent Then, there exists a pointed complete immersed submanifold $\Sigma_\infty=(S_\infty,q_\infty)$ of $M_\infty$ of type $C^{k-1,1}$
such that $i(q_\infty)=p_\infty$ and that, after extraction of a subsequence, $(\Sigma_n,q_n)_\ninn$ converges towards
$(\Sigma_\infty,q_\infty)$ in the pointed weak $C^{k-1,1}$ Cheeger/Gromov topology.
\endproclaim
\proclabel{TheoremTC}
\noindent The terms used in this theorem are explained in section $2$ and appendix \headref{HolderSpaces} of this paper. What we call the Cheeger/Gromov topology is essentially the canonical topology that one would expect to use for immersed submanifolds. In particular, when $k=2$, the condition on the submanifolds amounts to a bound on the norms of the second fundamental forms
of the immersed submanifolds.
\medskip
\noindent We would like to underline that in \cite{Corlette}, Corlette clearly states that he anticipates that his finiteness result should arise from a compactness result in a way analagous to the Cheeger/Gromov case. We consider however that, given the importance of this result to our own work, it was necessary to properly unearth this theorem and to state and prove it in its own right.
\medskip
\noindent In the second section, we introduce various concepts associated to immersed submanifolds, and we describe the Cheeger/Gromov topology for pointed Riemannian manifolds and immersed submanifolds. In the third section, we study the manner in which immersed submanifolds may locally be described in terms of graphs over the tangent space to each point. The objective here being to bound from below the radius of the disk over which the submanifold is a graph, and to bound from above the derivatives of the function of which the submanifold is a graph. In the fourth section, using the technical results of the third section, we prove theorem \procref{TheoremTC} and we then prove Corlette's result \cite{Corlette} as a corollary. 
\medskip
\noindent In this paper, we will be working in the $C^{k,\alpha}$ category. Since we geometers are
not in general in the habit of using mappings of this type, appendix \headref{HolderSpaces} reviews the properties 
required of a class of functions for one to be able to construct a theory of manifolds out of it, and we then show 
how the class of $C^{k,\alpha}_\oploc$ mappings satisfies these properties. Finally, in appendix B, we provide a proof
of the now classical compactness theorem of Riemannian geometry in a form that is most appropriate for our uses. We also hope (perhaps vainly) that we have provided here a slightly more accessible proof of what is an important, but technically challenging result.
\medskip
\noindent I would like to thank Fran\c{c}ois Labourie for introducing me to the subject and drawing my attention to the utility of this result.
\newhead{Convergence of Manifolds}
\newsubhead{Immersed Submanifolds}
\noindent We review the basic definitions from the theory of immersed submanifolds and establish the notations that will be used throughout this article.
\medskip
\noindent Let $M$ be a smooth manifold. An {\emph immersed submanifold\/} is a pair $\Sigma=(S,i)$ where $S$ is a smooth manifold and $i:S\rightarrow M$ is a smooth immersion. Let $g$ be a Riemannian metric on $M$. We give $S$ the unique Riemannian metric $i^*g$ which makes $i$ into an isometry. We say that $\Sigma$ is {\emph complete\/} if and only if the Riemannian manifold $(S,i^*g)$ is.
\medskip
\noindent We introduce the following definition:
\proclaim{Definition \procref{DefnHigherDerivativesSecondFF}}
\noindent Let $(M,g)$ be a Riemannian manifold. Let $X=(Y,i)$ be an immersed submanifold in $M$. Let $\nabla^i$ be the Levi-Civita covariant derivative generated over $Y$ by
the immersion $i$ into $(M,g)$. Let $A(X)$ be the second fundamental form of $X$. For all $k\geqslant 2$, we define $A_k(X)$ using the following recurrence relation:
$$\matrix
A_2(X) \hfill&= A(X), \hfill\cr
A_k(X) \hfill&= \nabla^i A_{k-1}(X)\ \forall k\geqslant 3.\hfill\cr
\endmatrix$$
\noindent We now define $\Cal{A}_k(X)$ for all $k\geqslant 2$ by:
$$
\Cal{A}_k(X) = \sum_{i=2}^k \|A_i(X)\|.
$$
\endproclaim
\newsubhead{The Cheeger/Gromov Topology}
\noindent A {\emph pointed Riemannian manifold\/} is a pair $(M,p)$ where $M$ is a Riemannnian manifold and $p$ is a point in $M$. If $(M,p)$ and $(M',p')$ are pointed
manifolds then a {\emph morphism\/} (or {\emph mapping\/}) from $(M,p)$ to $(M',p')$ is a (not necessarily even continuous) function from $M$ to $M'$ which sends $p$ to $p'$ and
is $C^\infty$ in a neighbourhood of $p$. The family of pointed manifolds along with these morphisms forms a category. In this section, we will discuss a notion of convergence
for this family. It should be borne in mind that even though this family is not a set, we may consider it as such. Indeed, since every manifold may be plunged into an infinite
dimensional real vector space, we may discuss, instead, the equivalent family of pointed finite dimensional submanifolds of this vector space, and this is a set.
\medskip
\noindent Let $(M_n,p_n)_{n\in\Bbb{N}}$ be a sequence of complete pointed Riemannian manifolds. For all $n$, we denote by $g_n$ the Riemannian metric over $M_n$. We say that the
sequence $(M_n,p_n)_{n\in\Bbb{N}}$ {\emph converges\/} to the complete pointed manifold $(M_0,p_0)$ in the {\emph Cheeger/Gromov topology\/} if and only if:
\medskip
\myitem{(i)} for all $n$, there exists a mapping $\varphi_n:(M_0,p_0)\rightarrow (M_n,p_n)$,
\medskip
\noindent such that, for every compact subset $K$ of $M_0$, there exists $N\in\Bbb{N}$ such that for all $n\geqslant N$:
\medskip
\myitem{(i)} the restriction of $\varphi_n$ to $K$ is a $C^\infty$ diffeomorphism onto its image, and
\medskip
\myitem{(ii)} if we denote by $g_0$ the Riemannian metric over $M_0$, then the sequence of metrics $(\varphi_n^*g_n)_{n\geqslant N}$ converges to $g_0$ in the
$C^\infty$ topology over $K$.
\medskip
\noindent We refer to the sequence $(\varphi_n)_{n\in\Bbb{N}}$ as a sequence of {\emph convergence mappings\/} of the sequence $(M_n, p_n)_{n\in\Bbb{N}}$ with respect to the
limit $(M_0,p_0)$. The convergence mappings are trivially not unique. However, two sequences of convergence mappings $(\varphi_n)_{n\in\Bbb{N}}$ and $(\varphi'_n)_{n\in\Bbb{N}}$
are equivalent in the sense that there exists an isometry $\phi$ of $(M_0,p_0)$ such that, for every compact subset $K$ of $M_0$, there exists $N\in\Bbb{N}$ such that for
$n\geqslant N$:
\medskip
\myitem{(i)} the mapping $(\varphi_n^{-1}\circ\varphi_n')$ is well defined over $K$, and
\medskip
\myitem{(ii)} the sequence $(\varphi_n^{-1}\circ\varphi_n')_{n\geqslant N}$ converges to $\phi$ in the $C^\infty$ topology over $K$.
\medskip
\noindent One may verify that this mode of convergence does indeed arise from a topological structure over the space of complete pointed manifolds. Moreover, this topology
is Haussdorf (up to isometries).
\medskip
\noindent Most topological properties are unstable under this limiting process. For example, the limit of a sequence of simply connected manifolds is not necessarily
simply connected. On the other hand, the limit of a sequence of surfaces of genus $k$ is a surface of genus at most $k$ (but quite possibly with many holes).
\medskip
\noindent Let $M$ be a complete Riemannian manifold. A {\emph pointed immersed submanifold\/} in  $M$ is a pair $(\Sigma,p)$ where $\Sigma=(S,i)$ is an immersed submanifold in
$M$ and $p$ is a point in $S$.
\medskip
\noindent Let $(\Sigma_n,p_n)_{n\in\Bbb{N}}=(S_n,p_n,i_n)_{n\in\Bbb{N}}$ be a sequence of complete pointed immersed submanifolds in $M$. We say that
$(\Sigma_n,p_n)_{n\in\Bbb{N}}$ {\emph converges\/} to $(\Sigma_0,p_0)=(S_0,p_0,i_0)$ in the {\emph Cheeger/Gromov topology\/} if and only if $(S_n,p_n)_{n\in\Bbb{N}}$ converges
to $(S_0,p_0)$ in the Cheeger/Gromov topology, and, for every sequence $(\varphi_n)_{n\in\Bbb{N}}$ of convergence mappings of $(S_n,p_n)_\ninn$ with respect to this limit, and
for every compact subset $K$ of $S_0$, the sequence of functions $(i_n\circ\varphi_n)_{n\geqslant N}$ converges to the function $(i_0\circ\varphi_0)$ in the $C^\infty$ topology
over $K$.
\medskip
\noindent As before, this mode of convergence arises from a topological structure over the space of complete immersed submanifolds. Moreover, this topology is
Haussdorf (up to isometries).
\medskip
\newsubhead{The Class of $C^{k,\alpha}$ Pointed Manifolds}
\noindent We define the space of $C^{k,\alpha}_\oploc$ mappings as in appendix \headref{HolderSpaces}, and we define the
following topological structure over the space $C^{k,\alpha}$:
\proclaim{Definition \nextprocno}
\noindent Let $\Omega\subseteq\Bbb{R}^n$ be an open set. Let $\suite{f}{n},f$ be functions over $\Omega$ of type $C^{k,\alpha}$. We say that $\suite{f}{n}$ converges to $f$ in the {\emph weak 
$C^{k,\alpha}$ topology} if and only if, for all $\beta\in(0,\alpha)$:
$$
(\|f_n-f\|_{C^{k,\beta}(\Omega)})_\ninn\rightarrow 0.
$$
\noindent If $(f_n)_\ninn,f$ are functions of type $C^{k,\alpha}_\oploc$ over $\Omega$, then we say that $\suite{f}{n}$ converges to $f$ in the {\emph weak $C^{k,\alpha}$ topology} if and only
if for all $p\in\Omega$ there exists a neighbourhood $V$ of $p$ in $\Omega$ such that
$\suite{f}{n}$ converges to $f$ in the weak $C^{k,\alpha}$ topology over $V$.
\endproclaim
\noindent We show in appendix \headref{HolderSpaces} that for $k\geqslant 1$, addition, multiplication, composition and
inversion of $C^{k,\alpha}_\oploc$ functions are continuous operations with respect to the weak $C^{k,\alpha}_\oploc$
topology.
\medskip
\noindent We show in appendix \headref{HolderSpaces} how the class of $C^{k,\alpha}_\oploc$ functions for $k\geqslant 1$
has sufficient structure for us to construct the class of $C^{k,\alpha}$ manifolds. We thus make the following definition:
\proclaim{Definition \nextprocno}
\noindent For $k\in\Bbb{N}$ and $\alpha\in (0,1]$, a $C^{k,\alpha}$ Riemannian manifold is a triplet $(S,\Cal{A},g)$
where:
\medskip
\myitem{(1)} $S$ is a connected topological manifold,
\medskip
\myitem{(2)} $\Cal{A}$ is a $C^{k,\alpha}$ atlas (i.e. all the transition mappings are of type $C^{k,\alpha}_\oploc$) and
\medskip
\myitem{(3)} $g$ is a $C^{k-1,\alpha}_\oploc$ metric over $(S,\Cal{A})$ (i.e. g is described
locally in every chart by a $C^{k-1,\alpha}_\oploc$ function).
\endproclaim
\noindent We define pointed $C^{k,\alpha}$ manifolds and immersed submanifolds as for smooth manifolds. We may also define
for such manifolds a topology analogous to the Cheeger/Gromov topology. In the case of convergence
of $C^{k,\alpha}$ manifolds, the convergence mappings are all of type $C^{k,\alpha}$ and the metrics converge in the weak $C^{k-1,\alpha}_\oploc$ topology. We call the resulting topology the {\emph weak $C^{k,\alpha}$ Cheeger/Gromov} topology.
\newsubhead{Uniqueness of The Cheeger/Gromov Limit}
\noindent Let $U\subseteq\Bbb{R}^n$ be an open subset of $\Bbb{R}^n$. Let $g$ be a Riemannian metric over $U$.
Let $TU$ be the tangent bundle over $U$. Let $\pi:TU\rightarrow U$ be the canonical projection. Let $TTU$ be the tangent
bundle over $TU$. Let $VTU$ be the vertical subbundle of $TTU$:
$$
VTU = \opKer(\pi).
$$
\noindent Let $HTU$ be the horizontal subbundle of $TTU$ associated to the Levi-Civita connection of $g$. We have:
$$
HTU \oplus VTU = TTU.
$$
\noindent We know that:
$$
HTU, VTU\cong \pi^*TU.
$$
\noindent Let us denote by $i_H$ (resp. $i_V$) the canonical isomorphism which sends $HTU$ (resp. $VTU$) into
$\pi^*TU$. We define the metric $Tg$ over $TTU$ such that:
$$\matrix
Tg|_{HTU}\hfill&= i_H^*\pi^*g,\hfill\cr
Tg|_{VTU}\hfill&= i_V^*\pi^*g,\hfill\cr
HTU \hfill&\perp_{Tg} VTU.\hfill\cr
\endmatrix$$
\noindent We call $Tg$ the {\emph Levi-Civita lifting} of $g$ over $TU$. Since the Levi-Civita connection of $g$ depends on the first
derivative of $g$, a $C^{k,\alpha}$ control on $g$ and $g^{-1}$ yields a $C^{k-1,\alpha}$ control on $Tg$ and $Tg^{-1}$.
\medskip
\noindent We have the following elementary result:
\proclaim{Lemma \nextprocno}
\noindent Let $U$,$V$ be open subsets of $\Bbb{R}^n$. Let $\varphi:U\rightarrow V$ be a diffeomorphism. Let $g$ and $h$ be Riemannian
metrics over $U$ and $V$ respectively. Let $M$ and $N$ be the matrices representing $g$ and $h$ respectively with respect to Euclidean
metric over $\Bbb{R}^n$. There exists a function $C_0$ such that if $\varphi^*h=g$, then:
$$
\|D\varphi\|_{C^0(U)} \leqslant C_0(\|M\|_{C^0(U)},\|M^{-1}\|_{C^0(U)}, \|N\|_{C^0(U)}, \|N^{-1}\|_{C^0(U)}).
$$
\endproclaim
\proclabel{LemmaBDONGBDONDPHI}
\proof Since $\varphi^*h=g$, we have:
$$
D\varphi^tND\varphi = M.
$$
\noindent Since the mapping $A\mapsto A^2$ defines a diffeomorphism of the space of symmetric positive definite matrices onto itself, we may
write:
$$
(N^{1/2}D\varphi)^t(N^{1/2}D\varphi) = M.
$$
\noindent For any matrix $A$, we have:
$$
\|A^tA\| = \|A\|^2.
$$
\noindent Thus:
$$
\|N^{1/2}D\varphi\|^2 = \|M\|.
$$
\noindent Consequently:
$$
\|D\varphi\| = \|N^{-1/2}N^{1/2}D\varphi\| \leqslant \|N^{-1/2}\|\|M^{1/2}\|.
$$
\noindent The result now follows.\qed
\medskip
\noindent Using the Levi-Civita lifting of the metric, we may generalise this result to higher derivatives of $\varphi$:
\proclaim{Lemma \nextprocno}
\noindent Let $U$,$V$ be open subsets of $\Bbb{R}^n$. Let $\varphi:U\rightarrow V$ be a diffeomorphism. Let $g$ and $h$ be Riemannian
metrics over $U$ and $V$ respectively. Let $M$ and $N$ be the metrics representing $g$ and $h$ respectively with respect to the
Euclidean metric over $\Bbb{R}^n$. For all $k\in\Bbb{N}$, there exists a function $C_k$ such that if $\varphi^*h=g$, then:
$$
\|D^k\varphi\|_{C^0(U)} \leqslant C_k(\|M\|_{C^k(U)},\|M^{-1}\|_{C^k(U)}, \|N\|_{C^k(U)}, \|N^{-1}\|_{C^k(U)}).
$$
\endproclaim
\proclabel{LemmaBDONDKGBDONDKPHI}
\proof We proof this by induction. By lemma \procref{LemmaBDONGBDONDPHI}, the result is true when $k=0$. Suppose that the result is true
for all $k\leqslant m$. Let $T^m\varphi$ be the $m$'th jet of $\varphi$ sending $T^mU$ into $T^mV$. Let $T^mg$ and $T^mh$ be the
$m$-fold Levi-Civita liftings of $g$ and $h$ respectively. We have:
$$
(T^m\varphi)^*T^mh = T^mg.
$$
\noindent Let $M_m$ and $N_m$ be the matrices of $T^mg$ and $T^mh$ respectively with respect to the Euclidean metric. Since $T^mU$ (resp.
$T^mV$) is a bundle over $U$ (resp. $V$), for $R\in\Bbb{R}^+$, we may consider the subbundle $B_R^mU$ (resp. $B_R^mV$) of balls of
radius $R$ (with respect to the Euclidean metric) in $T^mU$ (resp. $T^mV$). Since the result is true for all $k\leqslant m$, there exists
$R$ which depends only on $\|M\|_{C^m(U)}$,$\|M^{-1}\|_{C^m(U)}$,$\|N\|_{C^m(U)}$, $\|N^{-1}\|_{C^m(U)}$ such that:
$$
T^m\varphi(B_1^m(U)) \subseteq B_R^m(V).
$$
\noindent There exist functions $\hat{C}_{m,R}$ and $\hat{c}_{m,R}$ such that:
$$\matrix
\|N_m^{-1}\|_{C^0(B_R(V))} \hfill&\geqslant \hat{c}_{m,R}(\|N\|_{C^m(V)}, \|N^{-1}\|_{C^m(V)}),\hfill\cr
\|M_m\|_{C^0(B_1(U))} \hfill&\leqslant \hat{C}_{m,R}(\|M\|_{C^m(U)}, \|M^{-1}\|_{C^m(U)}).\hfill\cr
\endmatrix$$
\noindent Finally, there is a function $D_{m+1}$ such that:
$$
\|D^{m+1}\varphi\|_{C^0(U)} \leqslant D_{m+1}(\|DT^{m+1}\varphi\|_{C^0(B_1(U))}).
$$
\noindent The result now follows for $k=m+1$ by lemma \procref{LemmaBDONGBDONDPHI}, and the result follows for all $k$ by induction.\qed
\medskip
\goodbreak
\noindent In particular, we obtain the following corollary:
\proclaim{Corollary \nextprocno}
\noindent Let $U$, $V$ be open subsets of $\Bbb{R}^n$. Let $(g_n)_\ninn$,$g$ be Riemannian metrics over $U$ such that:
$$
(\|g_n-g\|_{C^k(U)})_\ninn\rightarrow 0.
$$
\noindent Let $(h_n)_\ninn$,$h$ be Riemannian metrics over $V$ such that:
$$
(\|h_n-h\|_{C^k(V)})_\ninn\rightarrow 0.
$$
\noindent Let $g_{\opEuc}$ be the Euclidean metric over $\Bbb{R}^n$ and suppose that there exists $\Lambda\in\Bbb{R}^+$ such that, for
all $n$:
$$
\frac{1}{\Lambda}g_\opEuc \leqslant g_n,h_n\leqslant \Lambda g_\opEuc.
$$
\noindent Let $(\varphi_n)_\ninn:U\rightarrow V$ be $C^{k+1}$ mappings such that, for all $n$:
$$
\varphi_n^*h_n = g_n,\qquad\varphi_n(0)=0.
$$
\noindent Then, there exists a $C^{k,1}$ mapping $\varphi_0:U\rightarrow V$ such that, after extraction of a subsequence $(\varphi_n)_\ninn$
converges to $\varphi_0$ in the weak $C^{k,1}$ topology. Moreover:
$$
\varphi_0^*h=g,\qquad\varphi_0(0)=0.
$$
\endproclaim
\proclabel{CorUniquenessOfCheegerGromov}
\proof This follows from the preceeding lemma and the classical Arzela-Ascoli theorem.\qed
\medskip
\remark This result yields the uniqueness of Cheeger/Gromov limits.
\newhead{Immersed Submanifolds Locally as Graphs}
\newsubhead{Immersed Submanifolds as Graphs}
\noindent It is trivial that every immersed submanifold may be described everywhere locally as a graph over a ball of a given radius in the tangent space to the submanifold at each point. In
this section, we will show how to obtain a bound from below for the radius of such a ball in terms of the norm of the second fundamental form in a neighbourhood of the given point. We then show how bounds on the derivatives of the function of which the submanifold is a graph may be obtained in terms of bounds on the derivatives of the second fundamental form of the
submanifold.
\medskip
\noindent In this section we will only consider $C^\infty$ manifolds, although the same reasoning remains valid for $C^{k,\alpha}$ manifolds. Let $(S,i)$ be an immersed submanifold in
$\Bbb{R}^n=\Bbb{R}^m\oplus\Bbb{R}^{n-m}$. Let $\Omega$ be an open subset of $S$ and let $V$ be an open subset of
$\Bbb{R}^m$. We say that $(\Omega,i)$ is a {\emph graph} over $V$ if and only if there exists a diffeomorphism
$\alpha:V\rightarrow\Omega$ and a function $f:V\rightarrow\Bbb{R}^{n-m}$ such that for all $x\in V$:
$$
i\circ\alpha(x) = (x,f(x)).
$$
\noindent We call $f$ the {\emph graph function} of $\Sigma$ and we call $\alpha$ the {\emph graph reparametrisation}
of $\Omega$.
\medskip
\noindent We have the following result:
\proclaim{Lemma \nextprocno}
\noindent Let $\Sigma=(S,i)$ be an immersed submanifold of $\Bbb{R}^n=\Bbb{R}^m\oplus\Bbb{R}^{n-m}$. Let $U_1,U_2$ be open subsets of $\Bbb{R}^m$. Suppose that there exist open subsets
$\Omega_1,\Omega_2$ of $S$ such that, for each $k$
$(\Omega_k,i)$ is a graph over $U_k$. For each $k$, let $\alpha_k:U_k\rightarrow\Omega_i$ be the graph reparametrisation of
$\Omega_k$ and let $f_k:U_i\rightarrow\Bbb{R}^{n-m}$ be the graph function.
\medskip
\noindent Suppose that there exists $p\in U_1\minter U_2$ such that:
$$
\alpha_1(p) = \alpha_2(p),
$$
\noindent then, for every $q$ in the connected component of $U_1\minter U_2$ containing $p$:
$$
\alpha_1(q) = \alpha_2(q),\qquad f_1(q) = f_2(q).
$$
\endproclaim
\proclabel{LemmaUniquenessOfGraphFunction}
\remark In other words, for a given pair of open sets $(U,\Omega)$, the graph function and the graph reparametrisation
are locally unique.
\medskip
\proof Let $V$ be the connected component of $U_1\minter U_2$ containing $p$. Let us define $X$ by:
$$
X = \left\{q\in V\text{ s.t. }\alpha_1(q) = \alpha_2(q)\right\}.
$$
\noindent The set $X$ is closed. Let $\pi:\Bbb{R}^n\rightarrow\Bbb{R}^m$ be the canonical projection. For all $q\in V$:
$$
(\pi\circ i)\circ\alpha_1(q) = (\pi\circ i)\circ\alpha_2(q) = q.
$$
\noindent Since $(\pi\circ i)$ is locally invertible, it follows that if $\alpha_1(q)=\alpha_2(q)$ then, for all
$q'$ in a neighbourhood $\Omega$ of $q$:
$$
\alpha_1(q) = \alpha_2(q).
$$
\noindent Consequently, $X$ is open. The result now follows.\qed
\medskip
\noindent We have the following definition:
\proclaim{Definition \nextprocno}
\noindent Let $(S,i)$ be an immersed submanifold in $\Bbb{R}^n$. Let $K$ be a closed subset of $S$, and let $\Omega$ be an open subset of
$\Bbb{R}^n$. We say that $(K,i)$ is {\emph complete with respect to $\Omega$} if and only if for every Cauchy sequence $\suite{x}{n}$ in
$S$ such that $(i(x_n))_\ninn$ converges in $\Omega$, there exists $x_0$ in $S$ such that $\suite{x}{n}$ converges to $x_0$.
\endproclaim
\remark If $\Omega$ is relatively compact, then this definition is independant of the Riemannian metric chosen over a neighbourhood of $\Omega$.
\medskip
\noindent Let $\Sigma=(S,i)$ be an immersed submanifold of $\Bbb{R}^n$. For $\epsilon\in\Bbb{R}^+$, suppose
that $\Sigma$ is a complete with respect to the ball $B_\epsilon(0)$. Let $p$ be a point in $S$ and suppose that:
$$
i(p) = 0,\qquad T_p\Sigma = i_*T_pS = \Bbb{R}^m\times\left\{0\right\}.
$$
\noindent We define $E$ to be the set of all $\eta\in (0,\infty)$ such that there exists a neighbourhood $\Omega$ of
$p$ such that $(\Omega,i)$ is a graph over $B_\eta(0)$. We define $\eta_0$ by:
$$
\eta_0 = \opSup(E).
$$
\noindent By the inverse function theorem, $E$ is non-empty, and consequently $\eta_0$ is well defined. Moreover, by lemma
\procref{LemmaUniquenessOfGraphFunction}, $\eta_0\in E$. We obtain the following result:
\proclaim{Lemma \nextprocno}
\noindent Let $r>\epsilon>0$ be positive real numbers, and suppose that the closed ball of radius $r$ about $p$ in $S$ is complete with respect to $B_\epsilon(0)$. There exists a function 
$\mu_{\text{graph}}(\epsilon,r)$ such that one of the following must be true:
\medskip
\myitem{(1)} $\eta_0\geqslant \epsilon/2$,
\medskip
\myitem{(2)} $\opSup\left\{\|f(p)\|\text{ s.t. }p\in B_{\eta_0}(0)\right\}\geqslant \epsilon/2$, or
\medskip
\myitem{(3)} $\opSup\left\{\|Df(p)\|\text{ s.t. }p\in B_{\eta_0}(0)\right\}\geqslant \mu_{\text{graph}}(\epsilon,r)$.
\endproclaim
\proclabel{LemmaDisasterConditions}
\remark There are three ways for the submanifold $\Sigma$ to stop being a graph:
\medskip
\myitem{(1)} The submanifold leaves $B_\epsilon(0)$.
\medskip
\myitem{(2)} The submanifold leaves $B_r(p)$.
\medskip
\myitem{(3)} The graph becomes vertical.
\medskip
\noindent The first two conditions in the lemma take into account the first form of degenation, wheras the last condition simultaneously takes into account the other two.
\medskip
\proof We define $\mu_{\text{graph}}(\epsilon,r)$ by:
$$
\mu_{\text{graph}}(\epsilon,r) = \sqrt{\frac{4r^2}{\epsilon^2}-1}.
$$
\noindent We will assume the contrary in order to obtain a contradiction. Let $\Omega$ be a neighbourhood of $p$, $\alpha:B_{\eta_0}(0)\rightarrow\Omega$ a diffeomorphism and $f:B_{\eta_0}(0)\rightarrow\Bbb{R}^{n-m}$ a mapping such that, for all
$q\in B_{\eta_0}(0)$:
$$
i\circ\alpha(q) = (q,f(q)).
$$
\noindent For all $q\in B_{\eta_0}(0)$, we have:
$$
\|Df(q)\| \leqslant \mu_{\text{graph}}(\epsilon,r).
$$
\noindent Since its derivative is bounded, $f$ extends to a continuous function on $\overline{B_{\eta_0}(0)}$. 
\medskip
\noindent Let us denote by $g_m$ the Euclidean metric over $\Bbb{R}^m$. Using the bound of $Df$, we find that there exists $\delta_2>0$ such
that:
$$
(i\circ\alpha)^* g_n \leqslant \left(\frac{4r^2}{\epsilon^2}-\delta_2\right)g_m.
$$
\noindent It follows that there exists $\delta_3>0$ such that for all $q\in B_{\eta_0}(0)$:
$$
\alpha(q) \in B_{r-\delta_3}(p).
$$
\noindent Moreover, if $\suite{q}{n}$ is a Cauchy sequence in $B_{\eta_0}(0)$ then $(\alpha(q_n))_\ninn$
is a Cauchy sequence in $\Omega$. Since the closed ball of radius $r$ about $p$ in $\Sigma$ is complete with respect to $B_\epsilon(0)$, it
follows that $\alpha$ extends to a continuous function over $\overline{B_{\eta_0}(0)}$.
\medskip
\noindent Now, there exists $\delta_1>0$ such that, for all $q\in\partial B_{\eta_0}(0)$:
$$\matrix
\|(i\circ\alpha)(q)\| \hfill&= \|q,f(q)\|\hfill\cr
&\leqslant \|q\| + \|f(q)\|\hfill\cr
&\leqslant\epsilon-\delta_1.\hfill\cr
\endmatrix$$
\noindent Consequently, for all $q\in\overline{B_{\eta_0}(0)}$:
$$
(q,f(q)) = (i\circ\alpha)(q)\in B_\epsilon(0).
$$
\noindent In otherwords, for all $q\in\overline{B_{\eta_0}(0)}$, $i(\alpha(q))$ is contained in $B_\epsilon(0)$.
\medskip
\noindent Let $\pi:\Bbb{R}^n\rightarrow\Bbb{R}^m\times\Bbb{R}^{n-m}$ be the canonical projection. Let $p$ be a point
in $\partial B_{\eta_0}(0)$. Since, for all $q\in\Omega$:
$$\matrix
\opDet(T_q(\pi\circ i)) \hfill&= \opDet(\delta_{ij}+\langle\partial_i f,\partial_j f\rangle\alpha^{-1}(q))^{-1/2}\hfill\cr
&\geqslant (1+\mu_{\text{graph}}(\epsilon,r)^2)^{-\frac{n}{2}}\hfill\cr
\endmatrix$$
\noindent It follows that:
$$
\opDet(T_{\alpha(p)}(\pi\circ i)) \neq 0.
$$
\noindent Thus, by the inverse function theorem, there exists $\eta_p\in\Bbb{R}^+$, a neighbourhood $\Omega_p$ of
$\alpha(p)$ in $S$, a diffeomorphism $\alpha:B_{\eta_p}(p)\rightarrow\Omega_p$ and a
function $f_p:B_{\eta_p}(p)\rightarrow\Bbb{R}^{n-m}$ such that for all $q$ in $B_{\eta_p}(p)$:
$$
(i\circ\alpha_p)(q) = (q,f_p(q)).
$$
\noindent Since $\alpha_p(p)\in\overline{\Omega}$, it follows that $\Omega\minter\Omega_p\neq\emptyset$. Let $q$ be a point
in $\Omega\minter\Omega_p$. We have:
$$
(\alpha\circ(\pi\circ i))(q) = (\alpha_p\circ(\pi\circ i))(q) = q.
$$
\noindent Since $(\pi\circ i)(q)\in B_{\eta_0}(0)\minter B_{\eta_p}(p)$ it follows by lemma
\procref{LemmaUniquenessOfGraphFunction} that $\alpha$ coincides with $\alpha_p$ over $B_{\eta_0}(0)\minter B_{\eta_p}(p)$.
\medskip
\noindent By compactness, there exists a finite set of points $p_1, ..., p_k\in\partial_{\eta_0}(p)$ and
$\delta\in\Bbb{R}^+$ such that:
$$
B_{\eta_0+\delta}(0) \subseteq B_{\eta_0}(0)\munion\left(\munion_{i=1}^kB_{\eta_{p_i}}(p_i)\right).
$$
\noindent Suppose that $B_{\eta_{p_i}(p_i)}\minter B_{\eta_{p_j}(p_j)}\neq\emptyset$. Since the straight line joining
$p_i$ and $p_j$, which is contained in $\overline{B_{\eta_0}(0)}$ intersects
$B_{\eta_{p_i}}(p_i)\minter B_{\eta_{p_j}}(p_j)$ non trivially, we have:
$$
B_{\eta_{p_i}}(p_i)\minter B_{\eta_{p_j}}(p_j)\minter \overline{B_{\eta_0}(0)} \neq\emptyset.
$$
\noindent The mappings $\alpha_{p_i}$ and $\alpha_{p_j}$ coincide in this set. Thus, since
$B_{\eta_{p_i}}(p_i)\minter B_{\eta_{p_j}}(p_j)$ is connected, it follows by lemma \procref{LemmaUniquenessOfGraphFunction}
that $\alpha_{p_i}$ and $\alpha_{p_j}$ coincide over this set.
\medskip
\noindent It thus follows that there exists $\delta\in\Bbb{R}^+$ such that we may extend $\alpha$ and $f$ to $C^{k,\alpha}$ mappings
defined on $B_{\eta_0+\delta}(0)$. Moreover, since:
$$
(\pi\circ i)\circ\alpha = \opId,
$$
\noindent it follows that $\alpha$ is a diffeomorphism onto its image. Consequently:
$$
\eta_0 + \delta\in E.
$$
\noindent We thus obtain a contradiction and the result now follows.\qed
\newsubhead{A Bound from Below of the Radius of Definition}
\noindent Let $\Omega$ be an open subset of $\Bbb{R}^n$. Let $g$ be a Riemannian metric over $\Omega$. Let $M:\Omega\rightarrow\opEnd(\Bbb{R}^n)$ be a smooth function taking values in the space of positive definite symmetric matrices
such that for all $V_p\in T_p\Omega$:
$$
g(V_p,V_p) = \langle V_p,M\cdot V_p\rangle.
$$
\noindent For every $p\in\Omega$, let $\opExp_p$ be the exponential mapping of $\Omega$ with respect to $g$ defined in a neighbourhood of $0\in T_p(\Omega)$. We recall the following facts concerning this application:
\proclaim{Lemma \nextprocno}
\noindent For all $p\in\Omega$:
$$\matrix
\|D\opExp_p(0)\| \hfill&= \|M^{-1}\|^{1/2},\hfill\cr
\|D\opExp_p(0)^{-1}\| \hfill&= \|M\|^{1/2}.\hfill\cr
\endmatrix$$
\noindent There exist functions $\mu_{\text{exp}}^1,\mu_{\text{exp}}^{-1}$ such that:
$$\matrix
\|D^2\opExp_p(0)\| \hfill&\leqslant \mu_{\text{exp}}^1(\|M\|,\|M^{-1}\|,\|DM\|),\hfill\cr
\|D^2(\opExp_p)^{-1}(p)\| \hfill&\leqslant \mu_{\text{exp}}^{-1}(\|M\|,\|M^{-1}\|,\|DM\|).\hfill\cr
\endmatrix$$
\endproclaim
\proclabel{LemmaControlOfExponential}
\noindent Next, we have the following result concerning the transformation of the second fundamental form:
\proclaim{Lemma \nextprocno}
\noindent There exists a continuous function $\mu_{II}^{\text{transform}}$ such that if:
\medskip
\myitem{(1)} $U$,$V$ are open sets,
\medskip
\myitem{(2)} $\varphi:U\rightarrow V$ is a diffeomorphism,
\medskip
\myitem{(3)} $\Sigma=(S,i)$ is an immersed hypersurface and,
\medskip
\myitem{(4)} $II$ (resp. $II'$) is the (Euclidean) second fundamental form of $(S,i)$ (resp. $(S,\varphi\circ i)$), 
\medskip
\noindent then, for all $p\in S$:
$$\matrix
\|II(p)\| \hfill&\leqslant \mu_{II}^{\text{transform}}(\|II'(p)\|,\|D\varphi(p)\|,\|D\varphi^{-1}(p)\|,\|D^2\varphi(p)\|),\hfill\cr
\|II'(p)\| \hfill&\leqslant \mu_{II}^{\text{transform}}(\|II(p)\|,\|D\varphi(p)\|,\|D\varphi^{-1}(p)\|,\|D^2\varphi(p)\|).\hfill\cr
\endmatrix$$
\noindent In particular, if $\|D\varphi\|$, $\|D\varphi^{-1}\|$ and $\|D^2\varphi\|$ are bounded, then a bound on
$\|II\|$ yields a bound on $\|II'\|$.
\endproclaim
\proclabel{LemmaTransformationsOfBoundsOnII}
\proof Since translations and rotations are isometries, we may assume that $\Sigma = (S,i)$
(resp. $\Sigma'=(S,\varphi\circ i)$) is the graph of a function $f$ (resp $f'$) such that $f(0)=df(0)=0$
(resp. $f'(0)=df'(0)=0$). In this case, $II$ (resp $II'$) coincides with $D^2f$ (resp $D^2f'$) and the result now
follows from a direct calculation.\qed
\medskip
\noindent By combining lemmata \procref{LemmaControlOfExponential} and \procref{LemmaTransformationsOfBoundsOnII} we now
obtain the following result:
\proclaim{Lemma \nextprocno}
\noindent There exists a continuous function $\mu_{II}^{\text{compare}}$ such that if $\Sigma$ is an immersed submanifold in $\Omega$ and $II^g$ (resp. $II$) is the second fundamental form of $\Sigma$ with respect to $g$ (resp. the Euclidean metric
over $\Omega$) then:
$$
\|II\| \leqslant \mu_{II}^{\text{compare}}(\|II^g\|,\|M\|,\|M^{-1}\|,\|DM\|).
$$
\noindent In particular, if $\|M\|$,$\|M^{-1}\|$ and $\|DM\|$ are bounded, then a bound on $II^g$ yields a bound on
$II$.
\endproclaim
\proclabel{LemmaBoundIIEuc}
\proof Let $p$ be a point in $S$. By applying $(\opExp_{i(p)})^{-1}$, we may work in an exponential chart about $i(p)$.
At the origin in such a chart, the second fundamental form of an immersed submanifold with respect to
$(\opExp_{i(p)})^*g$ coincides with the Euclidean second fundamental form. Using lemma
\procref{LemmaTransformationsOfBoundsOnII}, we may thus bound $\|II\|$ in terms of $\|II^g\|$ and the derivatives of $\opExp_{i(p)}$ at $p$. The result now follows by lemma 
\procref{LemmaControlOfExponential}.\qed
\medskip
\noindent We now obtain:
\proclaim{Lemma \nextprocno}
\noindent There exists a function $\mu_{II}(K,B,\epsilon,\|M\|,\|M^{-1}\|,\|DM\|)\leqslant\epsilon/2$
such that if:
\medskip
\myitem{(1)} $\delta\leqslant\mu_{II}$,
\medskip
\myitem{(2)} $B_\epsilon(0)\subseteq\Omega$,
\medskip
\noindent and if $f:B_\delta(0)\rightarrow\Bbb{R}^{n-m}$ is a function such that:
\medskip
\myitem{(1)} $f(0),df(0)=0$, and
\medskip
\myitem{(2)} the norm of the second fundamental form of the graph of $f$ with respect to $M$ is bounded above by $K$,
\medskip
\noindent then, for all $q\in B_\delta(0)$\colon
$$\matrix
\|f(q)\| \hfill&\leqslant \epsilon/2,\hfill\cr
\|Df(q)\| \hfill&\leqslant B.\hfill\cr
\endmatrix$$
\endproclaim
\proclabel{LemmaBoundOnDDf}
\remark In otherwords, we start by studying the graph of a function $f$ such that\break $f(0),df(0)=0$ 
(these conditions reflect the fact that, in the sequel, we will be studying immersed submanifolds in terms of graphs over the tangent space at each point). Then, given a bound $K$ on the second fundamental form of the graph, and given a desired bound $B$ on the derivative of $f$, we find a
radius $\delta$, depending only on $B$ and $K$ (and various other variables nonetheless 
independant of $f$), over which this bound is satisfied (provided, of course, that $M\circ f$ is defined
over this radius, hence $\delta<\epsilon/2$). Finally, for no extra cost, we also obtain a bound for $f$ over this radius which will be useful in the sequel.
\medskip
\proof By lemma \procref{LemmaBoundIIEuc}, the norm of the second fundamental form of the graph of $f$ with respect to
the Euclidean metric if bounded above by $K'=\mu_{II}^{\text{compare}}(K,\|M\|,\|M^{-1}\|,\|DM\|)$. For all $i$, we denote $\hat{\partial}_i$
and $\hat{\msf{N}}_i$ by:
$$
\hat{\partial}_i = \pmatrix \partial_i \cr Df\cdot\partial_i \cr\endpmatrix,\qquad
\hat{\msf{N}}_i = \pmatrix Df^t\cdot \partial_i \cr -\partial_i \cr\endpmatrix.\qquad
$$
\noindent $(\hat{\partial}_1, ...,\hat{\partial}_n)$ is a basis of tangent vectors to the graph of $f$ and
$(\hat{\msf{N}}_1, ...,\hat{\msf{N}}_{n-m})$ is a basis of normal vectors to the graph of $f$. Let $II$ be the second fundamental form
of the graph of $f$ with respect to the Euclidean metric. We have:
$$\matrix
\left|\langle II(\hat{\partial}_i,\hat{\partial}_j),\hat{\msf{N}}_k\rangle\right|\hfill&
=\left|\langle D_{\hat{\partial}_i}\hat{\msf{N}}_k, \hat{\partial}_j\rangle\right|\hfill\cr
&=\left|\partial_i\partial_jf^k\right|\hfill\cr
\endmatrix$$
\noindent Thus:
$$\matrix
&\left|\partial_i\partial_j f^k\right|\hfill&\leqslant K'\|\hat{\partial}_i\|\|\hat{\partial}_j\|\|\hat{\msf{N}}_k\|\hfill\cr
& &\leqslant K'(1 + \|Df\|^2)^{3/2} \hfill\cr
\Rightarrow\hfill& \|D^2f\|\hfill&\leqslant K'(1+\|Df\|^2)^{3/2}\hfill\cr
\endmatrix$$
\noindent Solving this differential inequality with the intial conditions $f(0),Df(0)=0$, we obtain the desired result.\qed
\medskip
\goodbreak
\noindent We now obtain the following result as an immediate corollary:
\proclaim{Lemma \nextprocno}
\noindent For $r>\epsilon>0$, there exist functions $\Delta(K,r,\epsilon,\|M\|,\|M^{-1}\|,\|DM\|)\leqslant\epsilon/2$, and\break
$B(r,\epsilon,\|M\|,\|M^{-1}\|)$ such that if
$(\Sigma,p)=(S,i,p)$ is a pointed immersed submanifold of $B_\epsilon(0)$ such that:
\medskip
\myitem{(1)} the closed ball of radius $r$ about $p$ in $\Sigma$ is complete with respect to $B_\epsilon(0)$,
\medskip
\myitem{(2)} $i(p)=0$, and
\medskip
\myitem{(3)} if $II$ is the second fundamental form of $\Sigma$ with respect to $g$, then:\par
$$
\|II\| \leqslant K,
$$
\noindent then, for all $\delta\leqslant\Delta$, there exists a unique neighbourhood $U$ of $p$ in $S$ such that $(U,i)$ is a graph over
a Euclidean ball of radius $\delta$ about the origin, and if $f$ is the graph function of $\Sigma$ over $B_\delta(0)$, then:
$$
\|f\|\leqslant\epsilon/2,\qquad \|Df\|\leqslant B.
$$
\endproclaim
\proclabel{LemmaGraphFunctionsOverGivenRadius}
\proof The ball of radius $\|M\|^{-1/2}r$ about $p$ in $S$ with respect to the Euclidean metric over $\Bbb{R}^n$ is contained within the
ball of radius $r$ about $p$ in $S$ with respect to the metric $g$. We thus choose:
$$
B = \mu_{\text{graph}}(\epsilon, \|M^{-1}\|^{1/2}r).
$$
\noindent We now define:
$$
\Delta = \mu_{II}(K,B,\epsilon,\|M\|,\|M^{-1}\|,\|DM\|).
$$
\noindent The result now follows from lemmata \procref{LemmaDisasterConditions} and \procref{LemmaBoundOnDDf}.\qed
\newsubhead{Bounds on the Higher Derivatives of the Graph Function}
\noindent In this section, we aim to show how bounds on the higher derivatives of a function may be obtained in terms of bounds on the higher derivatives
of the second fundamental form of its graph. The results of this section are essentially trivial, but we include them for completeness and
clarity.
\medskip
\noindent Let $\Omega$ be an open subset of $\Bbb{R}^n$. For every positive integer $k$, and for every $r\in\Bbb{R}^+$, we define
$B^k_r(0)$ to be the ball of radius $r$ about the origin in $\Bbb{R}^k$. Let $m$ be a positive integer not greater than $n$. Let
$\epsilon\in\Bbb{R}^+$ be a small positive real number and let $f:B^m_\epsilon(0)\rightarrow\Bbb{R}^{n-m}$ be a smooth function whose graph
is contained in $\Omega$. Let us denote the graph of $f$ by $\Sigma$.
\medskip
\noindent Let $\langle\cdot,\cdot\rangle$ be the Euclidean metric over $\Bbb{R}^n$. Let $g$ be a metric over $\Omega$ and let
$M:\Omega\rightarrow\opSymm(\Bbb{R}^n)$ be the matrix representing $g$ relative to the Euclidean metric. Thus, for all $V_p\in T\Omega$:
$$
g(V_p,V_p) = \langle V_p, M(p)\cdot V_p\rangle.
$$
\noindent Let $\partial_1,...,\partial_n$ be the canonical basis of $\Bbb{R}^n$. For all $1\leqslant i\leqslant m$, we define $\hat{\partial}_1,...,\hat{\partial}_m$ by:
$$
\hat{\partial}_i = (0,...,1,...,0,\partial_i f^{m+1},...,\partial_i f^{n})^t.
$$
\noindent For $m+1\leqslant i\leqslant n$, we define $\hat{\Cal{E}}_{m+1},...,\hat{\Cal{E}}_n$ by:
$$
\hat{\Cal{E}}_i = (-\partial_1 f^i, ..., -\partial_m f^i, 0, ...,1,...,0)^t.
$$
\noindent We define $\hat{\msf{N}}_{m+1}, ...,\hat{\msf{N}}_n$ by:
$$
\hat{\msf{N}}_i = M^{-1}\hat{\Cal{E}}_i.
$$
\noindent We find that $(\hat{\partial}_1, ...,\hat{\partial}_m)$ defines a moving (non-orthonormal) frame for $T\Sigma$. Similarly,
$(\hat{\msf{N}}_{m+1},...,\hat{\msf{N}_n})$ defines a moving (non-orthonormal) frame for the normal bundle to $\Sigma$ relative to metric
$g$. We define the matrix $B$ by:
$$
B = (\hat{\partial}_1, ...,\hat{\partial}_m, \hat{\msf{N}}_{m+1}, ...,\hat{\msf{N}}_n).
$$
\noindent Let $\nabla$ and $D$ be the Levi-Civita covariant derivatives of $g$ and the Euclidean metric respectively. Let $II$ be the second
fundamental form of $\Sigma$, the graph of $f$.
\medskip
\noindent In the sequel, for all $i,j,k$, we denote by $\Cal{F}(J^if, J^jM, J^kM^{-1})$ any function depending only on the derivatives of $f$, $M$ and $M^{-1}$ up to orders $i$,$j$ and $k$ respectively.
\medskip
\noindent We obtain the following result:
\proclaim{Lemma \nextprocno}
\noindent For all $i,j,k$:
$$
g(II(\hat{\partial}_i,\hat{\partial}_j),\hat{\msf{N}}_k)
=\partial_i\partial_jf^k + \Cal{F}(J^1f,J^1M,J^1M^{-1}).
$$
\endproclaim
\proclabel{LemmaControlZeroOrder}
\proof For all $i,j,k$, we have:
$$\matrix
g(II(\hat{\partial}_i,\hat{\partial}_j),\hat{\msf{N}}_k) \hfill&=
\langle \nabla_{\hat{\partial}_i}\hat{\partial}_j, M\hat{N}_k\rangle \hfill\cr
&=\langle D_{\hat{\partial}_i}\hat{\partial}_j,\hat{\Cal{E}}_k\rangle + \Cal{F}(J^1f,J^1M,J^1M^{-1})\hfill\cr
&=\partial_i\partial_jf +\Cal{F}(J^1f,J^1M,J^1M^{-1}).\hfill\cr
\endmatrix$$
\noindent The result now follows.\qed
\medskip
\noindent Let $\pi$ be the orthogonal projection onto $T\Sigma$ with respect to $g$. We have:
\proclaim{Lemma \nextprocno}
\noindent For all $V_p\in T\Sigma$:
$$
\pi(V_p) = \sum_{i=1}^n \langle (B^{-1})^t\partial_i,V_p\rangle \hat{\partial}_i.
$$
\endproclaim
\proclabel{LemmaProjOntoTGraph}
\proof Since $(B\partial_1, ..., B\partial_n)$ is a basis for $T_p\Omega$, there exists $a_1, ..., a_n\in\Bbb{R}$ such that:
$$
V_p = \sum_{i=1}^n a_i B\partial_i.
$$
\noindent Moreover, since $(\hat{\partial}_1, ...,\hat{\partial}_m)=(B\partial_1,...,B\partial_m)$ is a basis for $T\Sigma$ and since
$(\hat{\msf{N}}_{m+1}, ..., \hat{\msf{N}}_n)=(B\partial_{m+1}, ..., B\partial_n)$ is a basis for $T\Sigma^\perp$, we have:
$$\matrix
\pi(V_p) \hfill&= \sum_{i=1}^m a_i B\partial_i\hfill\cr
&= \sum_{i=1}^m a_i\hat{\partial}_i.\hfill\cr
\endmatrix$$
\noindent However, for $1\leqslant i\leqslant m$, we have:
$$\matrix
\langle (B^{-1})^t\partial_i, V_p\rangle \hfill&=
\langle \partial_i, B^{-1}\sum_{j=1}^na_jB\partial_j\rangle\hfill\cr
&=\langle\partial_i, \sum_{j=1}^n a_j\partial_j\rangle\hfill\cr
&=a_i.\hfill\cr
\endmatrix$$
\noindent The result now follows.\qed
\medskip
\noindent In particular, we obtain:
\proclaim{Corollary \nextprocno}
\noindent We have:
$$
\pi = \sum_{i=1}^m\Cal{F}(J^1f, J^0M, J^0M^{-1})\hat{\partial}_i.
$$
\endproclaim
\proclabel{CorSizeOfProjOntoTGraph}
\proof By definition of $B$:
$$
B = \Cal{F}(J^1f, J^0M, J^0M^{-1}).
$$
\noindent The result now follows by the preceeding lemma.\qed
\medskip
\noindent Let $\hat{\nabla}$ be the Levi-Civita covariant derivative of $\Sigma$ with respect to $g$. We obtain the following
generalisation of lemma \procref{LemmaControlZeroOrder}:
\proclaim{Lemma \nextprocno}
\noindent For all $a_1, ..., a_m,i,j,k$:
$$\matrix
g((\hat{\nabla}^mII)(\hat{\partial}_i,\hat{\partial}_j;\hat{\partial}_{a_1}, ...,\hat{\partial}_{a_m}),\hat{\msf{N}}_k)
\hfill&=\partial_{a_1}...\partial_{a_m}\partial_i\partial_jf^k\hfill\cr
&\qquad + \Cal{F}(J^{m+1}f, J^{m+1}M, J^{m+1}M^{-1})\hfill\cr.
\endmatrix$$
\noindent Moreover, for all $a_1, ..., a_m,i,j,k$:
$$
g((\hat{\nabla}^mII)(\hat{\partial}_i,\hat{\partial}_j;\hat{\partial}_{a_1}, ...,\hat{\partial}_{a_m}),\hat{\partial}_k)
=\Cal{F}(J^{m+1}f, J^{m+1}M, J^{m+1}M^{-1}).
$$
\endproclaim
\proclabel{LemmaControlAllOrders}
\proof We prove this result by induction. By lemma \procref{LemmaControlZeroOrder}, the result holds for $m=0$. For $m\geqslant 1$, we
have:
$$\matrix
g((\hat{\nabla}^mII)(\hat{\partial}_i,\hat{\partial}_j;\hat{\partial}_{a_1},...,\hat{\partial}_{a_m}),\hat{\msf{N}}_k)\hfill
&=g(\nabla_{\hat{\partial}_{a_m}}(\hat{\nabla}^{m-1}II)(\hat{\partial}_i,\hat{\partial}_j;\hat{\partial}_{a_1},...,\hat{\partial}_{a_{m-1}}),
\hat{\msf{N}}_k)\hfill\cr
&\qquad-g((\hat{\nabla}^{m-1}II)(\hat{\nabla}_{\hat{\partial}_{a_m}}\hat{\partial}_i,\hat{\partial}_j;
\hat{\partial}_{a_1},...,\hat{\partial}_{a_{m-1}}),\hat{\msf{N}}_k)\hfill\cr
&\qquad-...\hfill\cr
&=\partial_{a_m}g((\hat{\nabla}^{m-1}II)(\hat{\partial}_i,\hat{\partial}_j;\hat{\partial}_{a_1},...,\hat{\partial}_{a_{m-1}}),
\hat{\msf{N}}_k)\hfill\cr
&\qquad-g((\hat{\nabla}^{m-1}II)(\hat{\partial}_i,\hat{\partial}_j;\hat{\partial}_{a_1},...,\hat{\partial}_{a_{m-1}}),
\nabla_{\hat{\partial}_{a_m}}\hat{\msf{N}}_k)\hfill\cr
&\qquad-g((\hat{\nabla}^{m-1}II)(\hat{\nabla}_{\hat{\partial}_{a_m}}\hat{\partial}_i,\hat{\partial}_j;
\hat{\partial}_{a_1},...,\hat{\partial}_{a_{m-1}}),\hat{\msf{N}}_k)\hfill\cr
&\qquad-...\hfill\cr
\endmatrix$$
\noindent However:
$$
\nabla_{\hat{\partial}_{a_m}}\hat{\msf{N}}_k = \Cal{F}(J^2f, J^1M, J^1M^{-1}).
$$
\noindent Similarly:
$$\matrix
\hat{\nabla}_{\hat{\partial}_{a_m}}\hat{\partial}_i \hfill&= \pi(\nabla_{\hat{\partial}_{a_m}}\hat{\partial}_i)\hfill\cr
&=\pi(\Cal{F}(J^2f, J^1M, J^1M^{-1})).\hfill\cr
\endmatrix$$
\noindent Thus, by corollary \procref{CorSizeOfProjOntoTGraph}:
$$
\hat{\nabla}_{\hat{\partial}_{a_m}}\hat{\partial}_i =\sum_{j=1}^m\Cal{F}(J^2f,J^1M,J^1M^{-1})\hat{\partial}_i.
$$
\noindent Finally:
$$\matrix
\partial_{a_m}g((\hat{\nabla}^{m-1}II)(\hat{\partial}_i,\hat{\partial}_j;\hat{\partial}_{a_1}...\hat{\partial}_{a_{m-1}}),\hat{\msf{N}}_k) \hfill
&= \partial_{a_m}\Cal{F}(J^mf,J^mM,J^mM^{-1}) \hfill\cr
&= \Cal{F}(J^{m+1}f,J^{m+1}M,J^{m+1}M^{-1}).\hfill\cr
\endmatrix$$
\noindent The first result now follows by the induction hypothesis. The second result follows by a similar reasoning.\qed
\medskip
\noindent We now obtain the following result:
\proclaim{Lemma \nextprocno}
\noindent For every positive integer $m\geqslant 2$, there exists a function $B_m$ such that if $K\in\Bbb{R}^+$ satisfies:
$$
\|\Cal{A}^m_\Sigma\|\leqslant K,
$$
\noindent then:
$$
\|D^mf\|_{C^0(B_\epsilon(0))}\leqslant B_m(K,\|f\|_{C^1(B_\epsilon(0))}, \|M\|_{C^{m-1}(\Omega)}, \|M^{-1}\|_{C^{m-1}(\Omega)}).
$$
\endproclaim
\proclabel{LemmaControlInTermsOfII}
\proof Let $\|\cdot\|$ denote the Euclidean norm over $\Bbb{R}^n$. For all $i$, we have:
$$
\|\hat{\partial}_i\| \leqslant 1 + \|Df\|.
$$
\noindent Thus:
$$
g(\hat{\partial}_i,\hat{\partial}_i)^{1/2} \leqslant \|M\|^{1/2}(1+|\|Df\|).
$$
\noindent Similarly, for all $i$, we have:
$$
g(\hat{\msf{N}}_i,\hat{\msf{N}}_i)^{1/2}\leqslant \|M\|^{1/2}(1+\|Df\|).
$$
\noindent By lemma \procref{LemmaControlAllOrders}, we thus obtain, for all $a_1, ..., a_m$ and for all $k$:
$$\matrix
\left|\partial_{a_1}...\partial_{a_m}f^k\right|\hfill&\leqslant \left|g((\hat{\nabla}^{m-2}II)(\hat{\partial}_{a_1},...,\hat{\partial}_{a_m}),\hat{\msf{N}}_k)\right|\hfill\cr
&\qquad+\left|\Cal{F}(J^{m-1}f,J^{m-1}M,J^{m-1}M^{-1})\right|\hfill\cr
&\leqslant \|M\|^{(m+1)/2}(1 + \|Df\|)^{m+1}\Cal{A}^m_\Sigma + \left|\Cal{F}(J^{m-1}f,J^{m-1}M,J^{m-1}M^{-1})\right|\hfill\cr
\endmatrix$$
\noindent The result now follows by induction.\qed
\newhead{A Compactness Result for Immersed Submanifolds}
\newsubhead{Normalised Charts.}
\noindent We use the results of the preceeding section in order to obtain an Arzela-Ascoli type theorem which is the principal result of this paper. We begin with the following definition:
\proclaim{Definition \nextprocno}
\noindent Let $k\in\Bbb{N}$ be a positive integer and let $\rho,\Gamma\in(0,\infty)$ be positive
real numbers. Let $(M,g)$ be a Riemannian manifold of type at least $C^k$. For $p\in M$, a {\emph $(\Gamma,k)$-normalised chart of radius $\rho$} about $p$ is a $C^k$ coordinate chart $(\varphi,U,V)$ of $M$ such that:
\medskip
\myitem{(1)} $\varphi(p)=0$,
\medskip
\myitem{(2)} $B_\rho(0)$, the closed Euclidean ball of radius $\rho$ about $0$ in $\Bbb{R}^n$,  is contained in $V$, and
\medskip
\myitem{(3)} if $A$ is the matrix of the metric $\varphi_*g$ with respect to the Euclidean metric over $V$, then:\par
$$
\|A\|_{C^{k-1}(V)},\|A^{-1}\|_{C^{k-1}(V)} \leqslant \Gamma.
$$
\noindent For $p\in M$, we say that $M$ is {\emph $(\Gamma,k)$-normalisable over a radius $\rho$ about $p$} if and only if there
exists a $(\Gamma,k)$-normalised chart of $M$ of radius $\rho$ about $p$. For $K$ a subset of $M$, we say that $M$
is {\emph $(\Gamma,k)$-normalisable over a radius $\rho$ about $K$} if and only if $M$ is $(\Gamma,k)$-normalisable
over a radius $\rho$ about every point of $K$.
\endproclaim
\noindent We have the following elementary result:
\proclaim{Lemma \nextprocno}
\noindent For all $k\geqslant 1$, there exists a function $C_k$ such that if $(\varphi_i,U_i,V_i)_{i\in\left\{1,2\right\}}$ are
$(\Gamma,k)$-normalised charts, then, for every Euclidean ball $B$ contained in $\varphi_1(U_1\minter U_2)$, we have:
$$
\|\varphi_2\circ\varphi_1^{-1}\|_{C^k(B)} \leqslant C_k(\Gamma).
$$
\endproclaim
\proclabel{LemmaBDONTRANSITION}
\proof This follows directly from lemma \procref{LemmaBDONDKGBDONDKPHI}.\qed
\medskip
\noindent For all $k$, for all $r$, and for all $x\in\Bbb{R}^k$, we denote by $B^k_r(x)$ the Euclidean ball of radius $r$ about $x$ in
$\Bbb{R}^k$. Using the results of the preceeding sections, we show that, for any Riemannian manifold, $M$, given an immersed
submanifold $\Sigma=(S,i)$ of $M$, we may show that $\Sigma$ may be described everywhere locally as a graph of a function over
a disc in a normalised chart where the radius of the disc may be bounded below and where the $C^k$ norm of the function may be bounded from
above. Formally, we have the following technical lemma:
\proclaim{Lemma \nextprocno}
\noindent Let $k\geqslant 2, m\leqslant n\in\Bbb{N}$ be positive integers. Let
$\Gamma,B,\rho>0$ be positive real numbers. There exist positive real numbers $C,r>0$ such that if:
\medskip
\myitem{(1)} $(M,g)$ is a complete Riemannian manifold of dimension $n$ and of type at least $C^k$,
\medskip
\myitem{(2)} $K$ is a subset of $M$ about which $M$ is $(\Gamma,k)$-normalisable over a radius $\rho$, and
\medskip
\myitem{(3)} $\Sigma=(S,i)$ is a complete immersed submanifold of $M$ of dimension $m$ and of type $C^k$,
\medskip
\noindent and if $p$ is a point in $S$ such that:
\medskip
\myitem{(1)} $i(p)\in K$, and
\medskip
\myitem{(2)} $\left|\Cal{A}^k_\Sigma(q)\right|\leqslant B$ for all $q\in B_\rho(p)$,
\medskip
\noindent then, for every $(\Gamma,k)$-normalised chart $(\varphi,U,V)$ of $M$ of radius $\rho$ about $i(p)$,
there exists:
\medskip
\myitem{(1)} a function $f:B_r^n(0)\rightarrow\Bbb{R}^{m-n}$ of type $C^k$ such that $f(0),df(0)=0$,
\medskip
\myitem{(2)} an open set $\Omega$ of $S$ about $p$, and
\medskip
\myitem{(3)} a diffeomorphism $\alpha:(B^n_r(0),0)\rightarrow(\Omega,p)$,
\medskip
\noindent such that:
\medskip
\myitem{(1)} $(\Omega,i)$ is the graph of $f$ over $B^m_r(0)$ in the chart $(\varphi,U,V)$ with
graph diffeomorphism $\alpha$. In otherwords, there exists
a rotation $A$ of $\Bbb{R}^n$ such that, for all $x\in B_r^n(0)$:
$$
\varphi\circ\alpha(x) = A(x,f(x)),
$$
\noindent and,
\medskip
\myitem{(2)} $\|f\|_{C^k(B_r^m(0))}\leqslant C$.
\endproclaim
\proclabel{LemmaLA}
\remark This result amounts to a ``globalisation'' of the previous section. To be precise, in the previous
section, we obtained a local description of submanifolds of bounded second fundamental form, whereas, in this section, we transform this data into
an atlas of controlled charts.
\medskip
\remark It is very important to observe that the graph diffeomorphism $\alpha$ also allows us to
construct a normalised chart for $\Sigma$ about $p$.
\medskip
\proof Let $(\varphi,U,V)$ be a $(\Gamma,k)$-normalised chart of $M$ of radius $\rho$ about $p$. By composing $\varphi$ with an
isometry if necessary, we may suppose that:
$$
(T\varphi\circ Ti)\cdot T_pS = \Bbb{R}^m\oplus\left\{0\right\}.
$$
\noindent We choose $r$ such that:
$$
r < \Delta(B,\rho',\rho,\Gamma,\Gamma,\Gamma),
$$
\noindent as in lemma \procref{LemmaGraphFunctionsOverGivenRadius}. Since $\Sigma$ is complete, the existence of $f$ now follows. Moreover:
$$
\|Df\| \leqslant B(\rho',\rho,\Gamma,\Gamma).
$$
\noindent We now choose $C$ such that:
$$
C > B_m(K,B(\rho',\rho,\Gamma,\Gamma),\Gamma,\Gamma).
$$
\noindent The bound on the derivatives of $f$ now follows from lemma \procref{LemmaControlInTermsOfII}, and the result now follows.\qed
\medskip
\noindent We now have the following result:
\proclaim{Lemma \nextprocno}
\noindent Let $k\geqslant 1$ be a positive integer. Let $(M_n,p_n)_\ninn,(M_\infty,p_\infty)$ be complete pointed Riemannian manifolds of
type at least $C^k$, and suppose that $(M_n,p_n)_\ninn$ coverges towards $(M_\infty,p_\infty)$ in the pointed $C^k$ Cheeger/Gromov topology.
\medskip
\noindent For all $R>0$, there exists $D,r>0$ such that for all $n\in\Bbb{N}\munion\left\{\infty\right\}$, the manifold $M_n$ is
$(D,k)$-normalisable about a radius $r$ over $B_R(p_n)$.
\endproclaim
\proclabel{LemmaLB}
\proof By compactness of the closure of $B_{R+1}(p_\infty)$, there exists $D_1,r_1>0$ such that $M_\infty$ is $(\Delta_1,k)$-normalisable about a radius $r_1$ over
$B_{R+1}(p_\infty)$.
\medskip
\noindent Let $(\varphi_n)_\ninn$ be a sequence of $C^k$ convergence mappings of $(M_n,p_n)_\ninn$ with respect to
$(M_\infty,p_\infty)$. There exists $N\in\Bbb{N}$ such that for $n\geqslant N$:
\medskip
\myitem{(1)} the restriction of $\varphi_n$ to $B_{R+1}(p_\infty)$ is a diffeomorphism onto its image, and
\medskip
\myitem{(2)} $B_R(p_n)\subseteq \varphi_n(B_{R+1}(p_\infty))$.
\medskip
\noindent Since $(\varphi_n^*g_n)_\ninn\rightarrow g_\infty$ in the $C^k_\oploc$ topology, there exists $D_2\geqslant D_1$
and $r_2\leqslant r_1$ such that for all $n\geqslant N$, the manifold $M_n$ is $(D_2,k)$-normalisable about a radius
$r_2$ over $B_R(p_n)$. The result for $n<N$ follows trivially by compactness and the result now follows.\qed
\newsubhead{The Compactness Result}
\noindent We now obtain the following result:
\proclaim{Theorem \procref{TheoremTC}}
\noindent Let $k\geqslant 2, m\leqslant n\in\Bbb{N}$ be positive integers.
\medskip
\noindent Let $(M_n,p_n)_\ninn,(M_\infty,p_\infty)$ be complete pointed Riemannian manifolds of dimension $n$ and of class at least
$C^k$ such that $(M_n,p_n)_\ninn$ converges towards $(M_\infty,p_\infty)$ in the pointed $C^k$ Cheeger/Gromov topology.
\medskip
\noindent For all $n\in\Bbb{N}$, let $\Sigma_n=(S_n,q_n)$ be an $m$-dimensional pointed immersed submanifold of $M$ of class at least $C^k$ such that $i(q_n)=p_n$.
\medskip
\noindent Suppose that for all $R>0$, there exists $B$ such that, for all $n$:
$$
\left|\Cal{A}_{\Sigma_n}^k(q)\right|\leqslant B\qquad \forall q\in B_R(q_n).
$$
\noindent Then, there exists a pointed complete immersed submanifold $\Sigma_\infty=(S_\infty,q_\infty)$ of $M_\infty$ of type $C^{k-1,1}$
such that $i(q_\infty)=p_\infty$ and that, after extraction of a subsequence, $(\Sigma_n,q_n)_\ninn$ converges towards
$(\Sigma_\infty,q_\infty)$ in the pointed weak $C^{k-1,1}$ Cheeger/Gromov topology.
\endproclaim
\proof For all $n$, let $g_n$ be the metric over $M_n$ and let $g_\infty$ be the metric over $M_\infty$. By lemmata
\procref{LemmaLA} and \procref{LemmaLB}, for all $R\in\Bbb{R}$, there exists $K,D,D',r,r'\in\Bbb{R}^+$ and, for all $n$,
an atlas $\Cal{A}_{n,R}=(x_q,U_q,B_r(0))_{q\in B_R(q_n)}$ of $(D,k)$-normalised charts over a radius $r$ of $S_n$ 
and an atlas $\Cal{B}_{n,R}=(y_q,V_q,B_{r'}(0))_{q\in B_R(q_n)}$ of $(D',k)$-normalised charts over a radius $r'$ of
$M_n$ such that:
\medskip
\myitem{(1)} for every $n$ and for every $q$ in $B_R(q_n)$:
$$
x_q(q) = 0,
$$
\myitem{(2)} for every $n$ and for every $q$ in $B_R(q_n)$:
$$
y_q(i(q)) = 0,
$$
\myitem{(3)} for every $n$ and for every $q$ in $B_R(q_n)$, $(y_q\circ i_n\circ x_q^{-1})$ is defined over $B_r(0)$ and:
$$
\|y_q\circ i_n\circ x_q^{-1}\|_{C^k(B_r(q))} \leqslant K.
$$
\noindent Using lemma \procref{LemmaBDONTRANSITION}, we find that there exists $B\in\Bbb{R}^+$ such that, for every 
$n\in\Bbb{N}$, $\Cal{A}_{n,R}$ is a $(B,r)$-optimal $C^{k-1,1}$-atlas of $(S_n,q_n)$ over a radius $R$ (see definition
\procref{DefnOptimalAtlas}). Thus, by the compactness theorem of Riemannian geometry 
(theorem \procref{ThmConvergenceOfRiemannianGeometry}), there exists a complete pointed $C^{k-1,1}$-manifold
$(S_\infty,q_\infty)$ to which $(S_n,q_n)_\ninn$ converges, after extraction of a subsequence, in the weak $C^{k-1,1}$ 
Cheeger/Gromov topology.
\medskip
\noindent Let $\xi$ be a point in $S_\infty$. Let $(\xi_n)_\ninn\in(S_n)_\ninn$ be a sequence of points converging to $\xi$.
By compactness, we may assume that $i_n(\xi_n)$ converges to a point in $M_\infty$ that we will refer to, slightly abusively,
as $i_\infty(\xi_\infty)$. Let $(\varphi_n)_\ninn$ be a sequence of optimal $C^{k-1,1}$ convergence mappings of $(S_n,q_n)_\ninn$ with
respect to $(S_\infty,q_\infty)$. Let $(\psi_n)_\ninn$ be a sequence of convergence mappings of $(M_n,p_n)_\ninn$ with respect to
$(M_\infty,p_\infty)$. By definition of $(\varphi_n)_\ninn$, we may suppose that $(x_{\xi_n}\circ\varphi_n\circ x_\xi^{-1})_\ninn$
converges to the identity in the weak $C^{k-1,1}_\oploc$ topology. Similarly, we may suppose that 
$(y_{\xi_n}\circ\psi_n\circ y_{\xi}^{-1})_\ninn$ also converges to the identity in the weak $C^{k-1,1}_\oploc$ topology.
By the classical Arzela-Ascoli theorem, we may suppose that there exists a function $i_\infty:B_r(0)\rightarrow B_{r'}(0)$
of type $C^{k-1,1}_\oploc$ to which $(y_{\xi_n}\circ i_n\circ x_{\xi_n}^{-1})_\ninn$ converges in the weak
$C^{k-1,1}_\oploc$ topology. By composition and inversion, we thus find that 
$(y_\xi\circ\psi_n^{-1}\circ i_n\circ\varphi_n\circ x_\xi^{-1})_\ninn$ converges to $i_\infty$ in the weak 
$C^{k-1,1}_\oploc$ topology. By repeating this operation with $\xi$ taking values in a countable dense subset of $S_0$,
we find that there exists $i_\infty:S_\infty\rightarrow M_\infty$ to which $(i_n)_\ninn$ converges, after extraction of a subsequence, in 
the weak $C^{k-1,1}_\oploc$ topology.
\medskip
\noindent Let $\hat{g}_\infty$ and $g_\infty$ be the Riemannian metrics over $S_\infty$ and $M_\infty$ respectively. Since $(i_n^*g_n)_\ninn$ 
converges to $\hat{g}_\infty$ and to $i_\infty^*g_\infty$, we find that $i_\infty$ is an isometric immersion, and the result now follows.\qed
\medskip
\goodbreak
\noindent We obtain as a corollary to this result the following theorem of Corlette:
\proclaim{Theorem \nextprocno\ {\bf [Corlette, 1990]}}
\noindent Let $M$ be a compact manifold of dimension $n$. Let $B,K\in\Bbb{R}$ be two positive real numbers. Let $m\leqslant n\in\Bbb{N}$
be positive integers. For all $\alpha\in(0,1]$, there exist only finitely many $C^{1,\alpha}$ isotopy classes of compact immersed
submanifolds $\Sigma=(S,i)$ of $M$ such that:
$$\matrix
\|II_\Sigma\| \hfill&\leqslant B,\hfill\cr
\opDiam(\Sigma)\hfill&\leqslant K.\hfill\cr
\endmatrix$$
\endproclaim
\remark We observe that if $\|II_\Sigma\|\leqslant B$, then a bound on the volume of $\Sigma$ is equivalent to a bound on the diameter of
$\Sigma$. We thus obtain an analogous result for $\opVol(\Sigma)\leqslant K$.
\medskip
\proof We suppose the contrary and reason by absurdity. Let $(\Sigma_n,q_n)_\ninn=(S_n,i_n,q_n)$ be a sequence of pointed immersed
submanifolds of $M$ no two of which are $C^{1,\alpha}$ isotopy equivalent such that:
$$\matrix
\|II_\Sigma\| \hfill&\leqslant B,\hfill\cr
\opDiam(\Sigma) \hfill&\leqslant K.\hfill\cr
\endmatrix$$
\noindent By the compactness of $M$, after extraction of a subsequence, there exists $p_\infty$ such that
$(i_n(q_n))_\ninn$ converges to $p_\infty$. By theorem \procref{TheoremTC}, there exists a complete pointed immersed submanifold
$(\Sigma_\infty,p_\infty)$ of $M$ of type $C^{1,1}$ such that, after extraction of a subsequence, $(\Sigma_n,p_n)_\ninn$ converges
to $(\Sigma_\infty,p_\infty)$ in the weak $C^{1,1}_\oploc$ Cheeger/Gromov topology. To begin with, we have:
$$
\opDiam(\Sigma_\infty) \leqslant K.
$$
\noindent Consequently, $\Sigma_\infty$ is compact. It follows that there exists $N\in\Bbb{N}$ such that, for all $n\geqslant N$ there exists:
\medskip
\myitem{(1)} a $C^{1,\alpha}$ diffeomorphism $\varphi_n:S_\infty\rightarrow S_n$, and
\medskip
\myitem{(2)} an immersion $i'_n:S_\infty\rightarrow M$,
\medskip
\noindent such that:
\medskip
\myitem{(1)} for all $n\geqslant N$, $i_n'=i_n\circ\varphi_n$ and
\medskip
\myitem{(2)} $(i_n')_\ninn\rightarrow i_\infty$ in the weak $C^{1,\alpha}$ topology.
\medskip
\noindent Let $\epsilon$ be the radius of convergence of $M$. We may assume that for all $m,n\geqslant N$:
$$\matrix
&d(i_n'(p),i_m'(p)) \hfill&\leqslant\epsilon\qquad\forall p\in S_\infty\hfill\cr
\Rightarrow\hfill& d(i_n\circ(\varphi_n\circ\varphi_m^{-1})(p),i_m(p))\hfill&\leqslant\epsilon\qquad\forall p\in S_m\hfill\cr
\endmatrix$$
\noindent Consequently, if we denote by $N\Sigma_m$ the normal bundle over $\Sigma_m$ in $M$ and by $\opExp$ the exponential mapping of
$M_\infty$, we find that there exists a section $X$ of type $C^{1,\alpha}$ of $N\Sigma_m$ over $S_m$ such that:
$$
i_n\circ(\varphi_n\circ\varphi_m^{-1})(p) = \opExp(X(p))\qquad\forall p\in S_m.
$$
\noindent We define $(i_t)_{t\in[0,1]}$ by:
$$
i_t(p) = \opExp(tX(p))\qquad\forall p\in S_m.
$$
\noindent This defines a $C^{1,\alpha}$ isotopy between $i_m$ and $i_n\circ(\varphi_n\circ\varphi_m^{-1})$ and we thus obtain a contradiction.\qed
\inappendicestrue
\global\headno=0
\newhead{H\"older Spaces}
\newsubhead{Building Manifolds out of Classes of Functions}
\noindent The family of $C^\infty$ Riemannian manifolds is not sufficiently closed for our 
purposes. We are thus required to introduce the notion of $C^{k,\alpha}$ Riemannian manifolds,
where $k\geqslant 1$ and $\alpha$ is a real number in $(0,1]$. We begin by reviewing the properties required of a class
of functions for one to be able to construct a theory of manifolds out of it. Let $\Cal{C}$ be a functor which
associates to every pair of open sets $U\subseteq\Bbb{R}^m$ and $V\subseteq\Bbb{R}^n$ a topological space of functions
$\Cal{C}(U,V)$ sending $U$ into $V$ such that:
\headlabel{HolderSpaces}
\subheadlabel{ManifoldConditions}
\medskip
\myitem{(1)} $\Cal{C}$ is {\emph contravariant in $U$} and {\emph covariant in $V$}. Thus if $U'\subseteq U$ and 
if $V\subseteq V'$, then restriction defines a continuous mapping $\opRest:\Cal{C}(U,V)\rightarrow\Cal{C}(U',V')$.
\medskip
\myitem{(2)} $\Cal{C}(U,V)$ is contained in $C^1(U,V)$, the space of continuously differentiable mappings from $U$ to $V$. Moreover, this
inclusion is continous.
\medskip
\myitem{(3)} $\Cal{C}$ is {\emph locally defined}. In otherwords, if $\Cal{F}(U,V)$ denotes the family of functions 
sending $U$ into $V$, and if, for any open subset $U'$ of $U$ we denote:
$$
\Cal{F}(U,V)\minter\Cal{C}(U',V) = \left\{ f\in\Cal{F}(U,V)\text{ s.t. }f|_{U'}\in\Cal{C}(U',V)\right\},
$$
\noindent then, for any family $(U_a)_{a\in A}$ of open subsets of $U$, we obtain:
$$
U=\munion_{a\in A}U_a\Rightarrow \Cal{C}(U,V)=\minter_{a\in A}\Cal{F}(U,V)\minter\Cal{C}(U_a,V).
$$
\myitem{(4)} $\Cal{C}$ is {\emph reflexive}. In otherwords, if $f\in\Cal{C}(U,V)$ is a diffeomorphism, then $f^{-1}\in\Cal{C}(V,U)$. Moreover, this operation is {\emph continuous} in the sense that if $(f_n)_{n\in\Bbb{N}},f_0\in\Cal{C}(U,V)$, are diffeomorphisms such that $(f_n)_\ninn$ converges to $f_0$, then, for every relatively compact subset $W$ of $\opIm(f_0)$, there exists $N\in\Bbb{N}$ such that:
$$
n\geqslant N\Rightarrow W\subseteq\opIm(f_n),
$$
\noindent and $(f_n^{-1}|_W)_{n\geqslant N}$ converges to $(f_0^{-1})|_W$.
\medskip
\myitem{(5)} $\Cal{C}$ is {\emph transitive}. In otherwords, if $f\in\Cal{C}(U,V)$ and if $g\in\Cal{C}(V,W)$, then
$g\circ f\in\Cal{C}(U,W)$. Moreover, this operation is continuous. Thus, if $(f_n)_\ninn$ converges to $f_0$ and if
$(g_n)_\ninn$ converges to $g_0$, then $(g_n\circ f_n)_\ninn$ converges to $g_0\circ f_0$.
\medskip
\noindent Given such a functor $\Cal{C}$, we construct a class of manifolds, $\Cal{M}(\Cal{C})$, whose transition maps
are always in $\Cal{C}$.
\medskip
\myitem{(A)} Condition $(2)$ concerning differentiability could be replaced by continuity (one would then replace ``diffeomorphism'' by ``homeomorphism'' in condition $(4)$). A meaningful manifold theory can be constructed with simple 
continuity, but differentiability is essential if one wishes to introduce such tools as Riemannian metrics or
the implicit function theorem.
\medskip
\myitem{(B)} Since maps in $\Cal{C}$ are locally defined (condition $(3)$), we do not need to worry about the topology
of intersections of coordinate charts when we aim to show that a transition map is in $\Cal{C}$.
\medskip
\myitem{(C)} By reflexivity (condition $(4)$), it suffices to show that a transition map $\varphi$ is in $\Cal{C}(U,V)$
to know that its inverse $\varphi^{-1}$ is also in $\Cal{C}(U,V)$. This simplifies the construction of atlases.
\medskip
\myitem{(D)} By transitivity (condition $(5)$), in order to show that a given chart is compatible with a given atlas,
it suffices to show that it is compatible with a minimal family of charts in that atlas covering that chart. This also simplifies the construction of atlases.
\medskip
\myitem{(E)} Let $M$ and $N$ be two manifolds in $\Cal{M}(\Cal{C})$. Let $\varphi:M\rightarrow N$ be a continuous mapping.
Let $(U_a,x_a)_{a\in A}$ and $(V_b,y_b)_{b\in B}$ be atlases of $M$ and $N$ respectively such that:
$$
\forall a\in A\ \exists b(a)\in B\text{ s.t. }\varphi(U_a)\subseteq V_{b(a)}.
$$
\noindent By transitivity, if $(y_{b(a)}\circ\varphi\circ x_a^{-1})$ is in $\Cal{C}$ for all $a\in A$, then, for every pair
of charts $(U,x)$ and $(V,y)$ in $M$ and $N$ respectively such that $\varphi(U)\subseteq V$, the mapping
$y\circ\varphi\circ x^{-1}$ is also in $\Cal{C}$. This simplifies the construction of mappings of type $\Cal{C}$ between manifolds in
$\Cal{M}(\Cal{C})$. Moreover, transitivity allows us to show that the composition of two such mappings is also in $\Cal{C}$.
\medskip
\noindent Continuity of transitivity permits us to consider convergence of mappings between manifolds in $\Cal{M}(\Cal{C})$. Moreover, composition of mappings between manifolds in $\Cal{M}(\Cal{C})$ is continuous. By continuity of reflexivity, 
if $M,N\in\Cal{M}(\Cal{C})$ are of the same dimension, and if $(\varphi_n)_\ninn,\varphi_0:M\rightarrow N$ are 
homeomorphisms of type $\Cal{C}$ such that $(\varphi_n)_\ninn$ converges to $\varphi_0$, then $(\varphi_n^{-1})_\ninn$
converges to $\varphi_0^{-1}$. In otherwords, if $U$ is a relatively compact open subset of $\opIm(\varphi_0)$, then there
exists $N\in\Bbb{N}$ such that:
$$
n\geqslant N\Rightarrow U\in\opIm(\varphi_n),
$$
\noindent and $(\varphi_n^{-1}|_U)_\ninn$ converges to $(\varphi_0^{-1}|_U)$.
\medskip
\noindent This allows us to contruct a basic theory of manifolds. We now impose the following technical conditions
on $\Cal{C}$:
\medskip
\myitem{(6)} If $\Cal{A}(U,V)$ denotes the space of all affine mappings from $U$ to $V$ then:
$$
\Cal{A}(U,V)\subseteq\Cal{C}(U,V).
$$
\noindent Moreover, this inclusion is continuous.
\medskip
\myitem{(7)} $\Cal{C}$ is {\emph closed under Cartesian products}. In otherwords, if $f\in\Cal{C}(U,V)$ and 
$g\in\Cal{C}(U',V')$, then:
$$
f\times g\in\Cal{C}(U\times U',V\times V').
$$
\noindent Moreover, the inclusion $\Cal{C}(U,V)\times\Cal{C}(U',V')\rightarrow\Cal{C}(U\times U',V\times V')$ is continuous.
\medskip
\noindent These conditions allow us to use the implicit function theorem and Cartesian products to construct manifolds.
\medskip
\noindent The following conditions are also useful for technical local constructions:
\medskip
\myitem{(8)} $\Cal{C}(U,\Bbb{R})$ is an algebra over $\Bbb{R}$. Moreover, addition and multiplication are continuous
operations. This algebra contains automatically the multiplicative identity since $\Cal{C}(U,\Bbb{R})$ contains all 
affine mappings.
\medskip
\myitem{(9)} For all $n$, $\Cal{C}(U,\Bbb{R}^n)$ is a vector space over $\Bbb{R}$ and $\Cal{C}(U,\Bbb{R})$. Moreover
the action of $\Cal{C}(U,\Bbb{R})$ on $\Cal{C}(U,\Bbb{R}^n)$ is continuous.
\medskip
\noindent We also wish to study families of functions and tensors over manifolds in $\Cal{M}(\Cal{C})$. Let $\Cal{F}$ 
be a functor which associates to every open subset $U$ of $\Bbb{R}^n$ a {\emph topological algebra} of functions
from $U$ to $\Bbb{R}$ such that:
\medskip
\myitem{(1)} $\Cal{F}$ is {\emph contravariant}. In otherwords, if $U'\subseteq U$, then $\Cal{F}(U)$ is contained in
$\Cal{F}(U')$.
\medskip
\myitem{(2)} $\Cal{F}$ contains the constant mappings (and thus the multiplicative identity).
\medskip
\myitem{(3)} $\Cal{C}$ acts on $\Cal{F}$ by {\emph pull back}. In otherwords, if $\varphi\in\Cal{C}(U,V)$ and if
$f\in\Cal{F}(V)$, then $f\circ\varphi\in\Cal{F}(U)$.
\medskip
\myitem{(4)} The action of pull back is {\emph continuous}. In otherwords, if $(\varphi_n)_\ninn,\varphi_0\in\Cal{C}(U,V)$
are such that $(\varphi_n)_\ninn$ converges to $\varphi_0$, and if $(f_n)_\ninn,f_0\in\Cal{F}(V)$ are such that
$(f_n)_\ninn$ converges to $f_0$, then $(f_n\circ\varphi_n)_\ninn$ converges to $f_0\circ\varphi_0$.
\medskip
\noindent These conditions permit us to define, for a given $M\in\Cal{M}(\Cal{C})$, the topological algebra $\Cal{F}(M)$ of
real functions of class $\Cal{F}$ over $M$. Moreover, if $\varphi:M\rightarrow N$ is a mapping of type $\Cal{C}$, then pull back (composition) defines a continuous mapping from $\Cal{F}(N)$ to $\Cal{F}(M)$.
\medskip
\noindent Tensors over manifolds in $\Cal{M}(\Cal{C})$ are defined in a similar way.
\newsubhead{Lipschitz and Locally Lipschitz Mappings}
\noindent Let $\Omega$ be an open subset of $\Bbb{R}^n$. Let $f:\Omega\rightarrow\Bbb{R}^n$ be a
function defined over $\Omega$. For $\alpha\in(0,1]$, we define $\nslipa{f}$ by:
$$
\nslipa{f} = \msup_{x\neq y\in\Omega}\frac{\|f(x)-f(y)\|}{\|x-y\|^\alpha}.
$$
\noindent We say that $f$ is an {\emph $\alpha$-Lipschitz} mapping over $\Omega$ if and only if:
$$
\nslipa{f} < \infty.
$$
\noindent We say that $f$ is a {\emph locally $\alpha$-Lipschitz} mapping if and only if for every $p\in\Omega$, there
exists a neighbourhood $V$ of $p$ in $\Omega$ such that $f$ is an $\alpha$-Lipschitz mapping over $V$. We define
$\lipa(\Omega)$ (resp. $\lipa_\oploc(\Omega)$) to be the space of $\alpha$-Lipschitz (resp. locally $\alpha$-Lipschitz)
mappings over $\Omega$.
\medskip
\remark We recall the following properties concerning Lipschitz mappings:
\medskip
\myitem{(1)} Every $\alpha$-Lipschitz mapping is continuous. Similarly, every locally $\alpha$-Lipschitz mapping is continuous.
\medskip
\myitem{(2)} $\nslipa{f}=0$ if and only if $f$ is constant.
\medskip
\myitem{(3)} Let us define $\opDiam(\Omega)$ to be the diameter of $\Omega$:
$$
\opDiam(\Omega) = \msup_{x,y\in\Omega}\|x-y\|.
$$
\noindent For all $\beta\in(0,\alpha]$, we have:
$$
\|f\|_{\text{Lip}^\beta(\Omega)} \leqslant \nslipa{f}\opDiam(\Omega)^{\alpha-\beta}.
$$
\noindent Consequently, if $\opDiam(\Omega)<\infty$, then:
$$
\beta\leqslant\alpha \Rightarrow \opLip^\alpha(\Omega) \subseteq \opLip^\beta(\Omega).
$$
\noindent However, in general:
$$
\beta\leqslant\alpha \Rightarrow \opLip_\oploc^\alpha(\Omega) \subseteq \opLip_\oploc^\beta(\Omega).
$$
\myitem{(4)} We consider $\Omega$ with the Euclidean metric as a Riemannian manifold. Let $d$ be the distance function
on $\Omega$ generated by this Riemannian metric. We define $\opDil(\Omega)$ by:
$$
\opDil(\Omega) = \msup_{x\neq y\in\Omega}\frac{d(x,y)}{\|x-y\|}.
$$
\noindent If $f\in C^1(\Omega)$ then:
$$
\|f\|_{\opLip^1(\Omega)} \leqslant \|f\|_{C^1(\Omega)}\opDil(\Omega).
$$
\noindent Thus, if $\opDil(\Omega)<\infty$, then:
$$
C^1(\Omega) \subseteq \opLip^1(\Omega).
$$
\noindent However, in general:
$$
C^1(\Omega)\subseteq \opLip^1_\oploc(\Omega).
$$
\myitem{(5)} The two preceeding remarks show that the structure of $\lipa(\Omega)$ depends on the geometry of $\Omega$
wheras the structure of $\lipa_\oploc(\Omega)$ does not. In the sequel, we define:
$$
\Delta(\Omega) = \opMax(\opDil(\Omega),\opDiam(\Omega)).
$$
\noindent We following trivial lemma summarises the operations under which the space of Lipschitz mappings is closed:
\proclaim{Lemma \nextprocno}
\noindent We have the following composition rules:
\medskip
\myitem{(1)} If $f,g\in\lipa(\Omega)$ then $f+g\in\lipa(\Omega)$ and:
$$
\nslipa{f+g} \leqslant \nslipa{f} + \nslipa{g}.
$$
\myitem{(2)} If $f,g\in\lipa(\Omega)\minter C^0(\Omega)$ then $f\cdot g\in\lipa(\Omega)$ and:
$$
\nslipa{f\cdot g} = \nsczero{f}\nslipa{g} + \nslipa{f}\nsczero{g}.
$$
\myitem{(3)} Let $\Omega'$ be another open set. If $f\in\lipa(\Omega')$, if $g\in\lip{\beta}(\Omega)$ and if
$g(\Omega)\subseteq\Omega'$ then $f\circ g\in\lip{\alpha\beta}(\Omega)$ and:
$$
\|f\circ g\|_{\lip{\alpha\beta}(\Omega)} \leqslant \|f\|_{\lipa(\Omega')}\|g\|^\alpha_{\lip{\beta}(\Omega)}.
$$
\endproclaim
\proclabel{LemmaBasicLipschitzComposition}
\newsubhead{H\"older Spaces}
\noindent In this section, we introduce $C^{k,\alpha}$ and $C^{k,\alpha}_\oploc$ mappings for $k\geqslant 1$, and provide 
a brief review of their properties. In particular, we show how $C^{k,\alpha}_\oploc$ mappings satisfy all the conditions
specified in section \subheadref{ManifoldConditions}.
\medskip
\noindent We recall that $\Bbb{N}$ denotes the set of positive integers, and, in particular, that $0\notin\Bbb{N}$. For
$k\in\Bbb{N}$, for $\alpha\in(0,1]$ and for $f:\Omega\rightarrow\Bbb{R}^m$ a $C^k$ function, we
define $\nslka{f}$, the $\lka$ norm of $f$ by:
$$
\nslka{f} = \sum_{0\leqslant i\leqslant k}\nsczero{D^i f} + \nslipa{D^k f}.
$$
\noindent We observe that $\nslka{\cdot}$ does indeed define a norm over the space of $C^k$ functions over $\Omega$. We say that a function {\emph is of type $\lka$} if and only if it is $C^k$ and its $\lka$ norm is finite. Similarly, we say that a function {\emph is of type $\lkaloc$} if and only if it is $C^k$ and, for every $p\in\Omega$ there exists a
neighbourhood $V$ of $p$ in $\Omega$ over which $f$ is of type $\lka$. We denote by $\slka$ (resp. $\slkaloc$) the space of $\lka$ (resp. $\lkaloc$) functions over $\Omega$. $C^{k,\alpha}_\oploc$ trivially satisfies condition $(1)$. All functions in $C^{k,\alpha}_\oploc$ are differentiable, and so $C^{k,\alpha}_\oploc$ satisfies condition $(2)$. It is locally defined, and thus satisfies condition $(3)$. It contains all affine mappings and is closed under Cartesian
products, and thus satisfies conditions $(6)$ and $(7)$. It thus remains to show reflexivity (condition $(4)$), transitivity (condition $(5)$), the fact that $C^{k,\alpha}_\oploc(\Omega,\Bbb{R})$ is a topological algebra (condition $(8)$) and the fact that $C^{k,\alpha}_\oploc(\Omega,\Bbb{R}^n)$ is a topological vector space upon which $C^{k,\alpha}_\oploc(\Omega,\Bbb{R})$ acts continuously (condition $(9)$). Using lemma \procref{LemmaBasicLipschitzComposition}, we now
obtain the following trivial result:
\proclaim{Lemma \nextprocno}
\noindent For all $k\in\Bbb{N}$ there exist continuous functions $\nu_{\text{prod}}^{\text{bnd}}$,
$\nu_{\text{comp}}^{\text{bnd}}$ and $\nu_{\text{inv}}^{\text{bnd}}$ such that, for all $\Omega,\Omega'$ and for all 
$\alpha\in(0,1]$:
\medskip
\myitem{(1)} If $f,g\in\slka$ then $f+g\in\slka$ and:
$$
\nslka{f+g} \leqslant \nslka{f} + \nslka{g}.
$$
\myitem{(2)} If $f,g\in\slka$ then $f\cdot g\in\slka$ and:
$$
\nslka{f\cdot g} \leqslant \nu_{\text{prod}}^{\text{bnd}}(\nslka{f},\nslka{g},\alpha,\Delta(\Omega)).
$$
\myitem{(3)} If $f\in\lka(\Omega')$, if $g\in\slka$ and if $g(\Omega)\subseteq\Omega'$,
then $f\circ g\in\slka$ and:
$$
\nslka{f\circ g} \leqslant \nu_{\text{comp}}^{\text{bnd}}(\|f\|_{C^{k,\alpha}(\Omega)},\nslka{g},\alpha,\Delta(\Omega)).
$$
\myitem{(4)} If $f:\Omega\rightarrow\Omega'$ is a $C^{k,\alpha}$ diffeomorphism, then:
$$
\|f^{-1}\|_{C^{k,\alpha}(\Omega)} \leqslant \nu_{\text{inv}}^{\text{bnd}}(\|f\|_{C^{k,\alpha}(\Omega)},
\|f^{-1}\|_{C^1(\Omega)},\alpha,\Delta(\Omega),\Delta(\Omega')).
$$
\endproclaim
\proclabel{LemmaAdvancedLipschitzComposition}
\noindent This now yields the following trivial corollary:
\proclaim{Corollary \nextprocno}
\myitem{(1)} If $\Delta(\Omega)<\infty$, then $\slka$ is an algebra over $\Bbb{R}$.
\medskip
\myitem{(1')} $\slkaloc$ is an algebra over $\Bbb{R}$.
\medskip
\myitem{(2)} If $\Delta(\Omega),\Delta(\Omega')<\infty$ and if $\varphi:\Omega\rightarrow\Omega'$ is a $\lka$ mapping, then:
$$
\varphi^*\lka(\Omega') \subseteq \slka.
$$
\medskip
\myitem{(2')}\noindent If $\varphi:\Omega\rightarrow\Omega'$ is a $\lkaloc$ mapping, then:
$$
\varphi^*\lkaloc(\Omega') \subseteq \slkaloc.
$$
\myitem{(3)} If $\Delta(\Omega),\Delta(\Omega'),\Delta(\Omega'')<\infty$ and if $\varphi:\Omega\rightarrow\Omega'$ and $\psi:\Omega'\rightarrow\Omega''$ are $\lka$ mappings, then $\psi\circ\varphi$ is also a $C^{k,\alpha}$ mapping.
\medskip
\myitem{(3')} If $\varphi:\Omega\rightarrow\Omega'$ and $\psi:\Omega'\rightarrow\Omega''$ are $\lkaloc$ mappings,
then $\psi\circ\varphi$ is a also a $\lkaloc$ mapping.
\medskip
\myitem{(4)} If $\Delta(\Omega),\Delta(\Omega')<\infty$ and if $\varphi:\Omega\rightarrow\Omega'$ is a $C^{k,\alpha}$ diffeomorphism, and if $\|\varphi^{-1}\|_{C^1(\Omega)}<\infty$, then $\varphi^{-1}$ is a $C^{k,\alpha}$ diffeomorphism.
\medskip
\myitem{(4')} If $\varphi:\Omega\rightarrow\Omega'$ is a $C^{k,\alpha}_\oploc$ diffeomorphism, then $\varphi^{-1}$ is
a $C^{k,\alpha}_\oploc$ diffeomorphism.
\endproclaim
\remark In particular, we see that $C^{k,\alpha}_\oploc(\Omega,\Bbb{R})$ is an algebra and that $C^{k,\alpha}_\oploc(\Omega,\Bbb{R}^n)$ is a vector space upon which $C^{k,\alpha}_\oploc(\Omega,\Bbb{R})$ acts by multiplication. Moreover, we see that $C^{k,\alpha}_\oploc$ is reflexive and transitive. It remains to show, however, that $C^{k,\alpha}_\oploc(\Omega,\Bbb{R})$ is a topological algebra, that $C^{k,\alpha}_\oploc(\Omega,\Bbb{R}^n)$ is a topological vector space, and
that $C^{k,\alpha}_\oploc(\Omega,\Bbb{R})$ acts continuously on this vector space. We must also show that reflexivity and
transitivity are continuous. This will be done in the next section. 
\medskip
\remark Although we are more interested in studying $C^{k,\alpha}_\oploc$ functions, when we are dealing with lemmata
comparing norms, it is preferable to study the family of $\lka$ functions, since, in this case, the results are slightly simpler to state.
\newsubhead{Continuity of Addition, Multiplication, Composition and Inversion}
\noindent We now study the topologies of $\slka$ and $\slkaloc$. Using lemma \procref{LemmaBasicLipschitzComposition}, we obtain the following result:
\proclaim{Lemma \nextprocno}
\myitem{(1)} If $f,f',g,g'\in\slipa$, then:
$$
\nslipa{(f+g)-(f'+g')} \leqslant \nslipa{f-f'} + \nslipa{g-g'}.
$$
\myitem{(2)} If $f,f',g,g'\in\slipa\minter C^0(\Omega)$, then:
$$\matrix
\nslipa{f\cdot g - f'\cdot g'} \hfill&\leqslant \nsczero{f}\nslipa{g-g'} + \nslipa{f}\nsczero{g-g'}\hfill\cr
&\qquad +\nsczero{f - f'}\nslipa{g'} + \nslipa{f-f'}\nsczero{g'}.\hfill\cr
\endmatrix$$
\endproclaim
\proclabel{LemmaBasicCtyAddMult}
\noindent This yields the following result which shows us that addition and multiplication are
continuous in the $\lka$ (resp. $\lkaloc$) topology over $\slka$ (resp. $\slkaloc$).
\proclaim{Lemma \nextprocno}
\myitem{(1)} If $f,f',g,g'\in\slka$, then:
$$
\nslka{(f+g)-(f'+g')}\leqslant \nslka{f-f'} + \nslka{g-g'}.
$$
\myitem{(2)} There exists a continuous function $\nu^{\text{cty}}_{\text{prod}}$ such that:\par
$$\matrix
\nslka{f\cdot g -  f'\cdot g'} \hfill&\leqslant \nu^{\text{cty}}_{\text{prod}}(\nslka{f-f'},\nslka{g-g'},\hfill\cr
&\qquad\qquad\nslka{f},\nslka{f'},\nslka{g},\hfill\cr
&\qquad\qquad\nslka{g'},\alpha,\Delta(\Omega)).\hfill\cr
\endmatrix$$
\noindent Moreover, $\nu^{\text{cty}}_{\text{prod}}(0,0,\cdot,\cdot,\cdot,\cdot,\cdot,\cdot)=0$.
\endproclaim
\proclabel{LemmaCtyAddMult}
\noindent In order to define a topology with respect to which composition of functions is
continuous, we require the following definition:
\proclaim{Definition \nextprocno}
\noindent Let $\Omega\subseteq\Bbb{R}^n$ be an open set. Let $\suite{f}{n},f$ be functions over $\Omega$ of type $C^{k,\alpha}$. We say that $\suite{f}{n}$ converges to $f$ in the {\emph weak 
$C^{k,\alpha}$ topology} if and only if, for all $\beta\in(0,\alpha)$:
$$
(\|f_n-f\|_{C^{k,\beta}(\Omega)})_\ninn\rightarrow 0.
$$
\noindent If $(f_n)_\ninn,f$ are functions of type $C^{k,\alpha}_\oploc$ over $\Omega$, then we say that $\suite{f}{n}$ converges to $f$ in the {\emph weak $C^{k,\alpha}_\oploc$ topology} if and only
if for all $p\in\Omega$ there exists a neighbourhood $V$ of $p$ in $\Omega$ such that
$\suite{f}{n}$ converges to $f$ in the weak $C^{k,\alpha}$ topology over $V$.
\endproclaim
\noindent We immediately see that $C^{k,\alpha}_\oploc(\Omega,\Bbb{R})$ is a topological algebra with respect to the weak $C^{k,\alpha}_\oploc$ topology. Moreover, $C^{k,\alpha}_\oploc(\Omega,\Bbb{R}^n)$ is a topological vector space, and
$C^{k,\alpha}_\oploc(\Omega,\Bbb{R})$ acts continuously on this space by multiplication. We now have the following result:
\proclaim{Lemma \nextprocno}
\myitem{(1)} If $f\in\lipa(\Omega')$, if $g,g'\in\opLip^1(\Omega)\minter C^0(\Omega)$ and if $g(\Omega), g'(\Omega)\subseteq\Omega'$, then,
for all $\beta\in (0,\alpha)$:
$$
\|f\circ g - f\circ g'\|_{\opLip^\beta(\Omega)} \leqslant
2^{1-\frac{\beta}{\alpha}}\|f\|_{\opLip^\alpha(\Omega')}\|g-g'\|^{\alpha-\beta}_{C^0(\Omega)}
(\|g\|_{\opLip^1(\Omega)}^\alpha + \|g'\|_{\opLip^1(\Omega)}^\alpha)^\frac{\beta}{\alpha}.
$$
\myitem{(2)} If $f,f'\in\opLip^\alpha(\Omega')\minter C^0(\Omega')$, if $g\in\opLip^1(\Omega)$ and if $g(\Omega)\subseteq\Omega'$, then, for all $\beta\in (0,\alpha)$:
$$
\|f\circ g - f'\circ g\|_{\opLip^\beta(\Omega)} \leqslant
2^{1-\frac{\beta}{\alpha}}\|f-f'\|_{C^0(\Omega')}^{1-\frac{\beta}{\alpha}}
(\|f\|_{\opLip^\alpha(\Omega')} + \|f'\|_{\opLip^\alpha(\Omega')})^{\frac{\beta}{\alpha}}\|g\|_{\opLip^1(\Omega)}^\beta.
$$
\endproclaim
\proclabel{LemmaBasicCtyComposition}
\remark The preceeding lemma is a special case of the following more general result:
\proclaim{Lemma \nextprocno}
\myitem{(1)} If $f\in\opLip^\alpha(\Omega)$, if $g,g'\in\opLip^\beta(\Omega')\minter C^0(\Omega)$ and if
$g(\Omega), g'(\Omega)\subseteq\Omega'$, then, for all $\lambda\in [0,1]$:
$$
\|f\circ g - f\circ g'\|_{\opLip^{\alpha\beta\lambda}(\Omega)} \leqslant
2^{1-\lambda}\|f\|_{\opLip^\alpha(\Omega')}\|g-g'\|_{C^0(\Omega)}^{\alpha(1-\lambda)}
(\|g\|_{\opLip^\beta(\Omega)}^\alpha + \|g'\|_{\opLip^\beta(\Omega)}^\alpha)^\lambda.
$$
\myitem{(2)} If $f\in\opLip^\alpha(\Omega')\minter C^0(\Omega')$, if $g\in\opLip^\beta(\Omega)$ and if $g(\Omega)\subseteq\Omega'$ then
for all $\lambda\in[0,1]$:
$$
\|f\circ g - f'\circ g\|_{\opLip^{\alpha\beta\lambda}} = 
2^{1-\lambda}\|f-f'\|_{C^0(\Omega')}^{1-\lambda}
(\|f\|_{\opLip^\alpha(\Omega')} + \|f'\|_{\opLip^\alpha(\Omega')})^\lambda\|g\|_{\opLip^\beta(\Omega)}^{\lambda\alpha}.
$$
\endproclaim
\noindent We now have:
\proclaim{Lemma \nextprocno}
\myitem{(1)} There exists a continuous function $\nu^{\text{cty}}_{\text{comp},1}$ such that if $f\in C^{k,\alpha}(\Omega')$, if
$g,g'\in C^{k,\alpha}(\Omega)$, if $g(\Omega),g'(\Omega)\subseteq\Omega'$, and if $\beta<\alpha$, then:\par
$$\matrix
\|f\circ g - f\circ g'\|_{C^{k,\beta}(\Omega)} \hfill&\leqslant
\nu^{\text{cty}}_{\text{comp},1}(\|g-g'\|_{C^{k,\alpha}(\Omega)}, \|g\|_{C^{k,\alpha}(\Omega)}, \|g'\|_{C^{k,\alpha}(\Omega)},\hfill\cr
&\qquad\qquad  \|f\|_{C^{k,\alpha}(\Omega')}, \alpha, \beta, \Delta(\Omega), \Delta(\Omega')).\hfill\cr
\endmatrix$$
\noindent Moreover $\nu^{\text{cty}}_{\text{comp},1}(0,\cdot,\cdot,\cdot,\cdot, \cdot, \cdot, \cdot)=0$.
\medskip
\noindent Consequently, if $\suite{g}{n}$ converges to $g_0$ in the $C^{k,\alpha}$ (resp. $C^{k,\alpha}_\oploc$) topology, and if\break
$\Delta(\Omega),\Delta(\Omega')<\infty$, then $(f\circ g_n)_\ninn$ converges to $(f\circ g_0)$ in the $C^{k,\beta}$ (resp. $C^{k,\beta}_\oploc$) topology.
\medskip
\myitem{(2)} There exists a continuous function $\nu^{\text{cty}}_{\text{comp},2}$ such that, if $f,f'\in C^{k,\alpha}(\Omega')$, if
$g\in C^{k,\alpha}(\Omega)$, if $g(\Omega)\subseteq\Omega'$, and if $\beta<\alpha$, then:\par
$$\matrix
\|f\circ g - f'\circ g\|_{C^{k,\beta}(\Omega)} \hfill&\leqslant
\nu^{\text{cty}}_{\text{comp},2}(\|f-f'\|_{C^{k,\alpha}(\Omega')}, \|f\|_{C^{k,\alpha}(\Omega')}, \|f'\|_{C^{k,\alpha}(\Omega')},\hfill\cr
&\qquad\qquad \|g\|_{C^{k,\alpha}(\Omega)}, \alpha, \beta, \Delta(\Omega), \Delta(\Omega')).\hfill\cr
\endmatrix$$
\noindent Moreover, $\nu^{\text{cty}}_{\text{comp},2}(0,\cdot,\cdot,\cdot,\cdot,\cdot,\cdot,\cdot)=0$.
\medskip
\noindent Consequently, if $\suite{f}{n}$ converges to $f_0$ in the $C^{k,\alpha}$ (resp. $C^{k,\alpha}_\oploc$)
topology, and if\break $\Delta(\Omega),\Delta(\Omega')<\infty$, then $(f_n\circ g)_\ninn$ converges to $(f_0\circ g)$ in the $C^{k,\beta}$
(resp. $C^{k,\beta}_\oploc$) topology.
\endproclaim
\proclabel{LemmaCtsComp}
\noindent This result shows us that the transitivity of $C^{k,\alpha}_\oploc$ is continuous. We observe that in order for transitivity to be continuous, we must have $k\geqslant 1$.
\medskip
\noindent Finally, using the chain rule, and the result for composition, we may also obtain the following result concerning inverses:
\proclaim{Lemma \nextprocno}
\noindent Let $\Omega,\Omega'\subseteq\Bbb{R}^n$ be open sets. Let $(\varphi_n)_\ninn,\varphi_0:\Omega\rightarrow\Omega'$ be $C^{k,\alpha}$ mappings which are diffeomorphisms onto their images. Suppose that $(\varphi_n)_\ninn$ converges to
$\varphi_0$ in the weak $C^{k,\alpha}_\oploc$ topology. Then, for every compact subset, $K$, of $\opIm(\varphi_0)$, there
exists $N\in\Bbb{N}$ such that:
$$
n\geqslant N\Rightarrow K\subseteq\opIm(\varphi_n),
$$
\noindent and $(\varphi^{-1}_n)_{n\geqslant N}$ converges to $\varphi_0^{-1}$ in the weak $C^{k,\alpha}_\oploc$ topology over $K$.
\endproclaim
\noindent It thus follows that reflexivity of $C^{k,\alpha}_\oploc$ is continuous. Again, in order for this to work, we require that $k\geqslant 1$.
\medskip
\noindent We have thus shown that for $k\geqslant 1$, $C^{k,\alpha}_\oploc$ functions satisfy all the conditions required for us to be able to construct a theory of $C^{k,\alpha}_\oploc$ manifolds. Similarly, the same results permit us to construct spaces of $C^{l,\beta}$ functions over such manifolds for $l+\beta\leqslant k+\alpha$ and spaces of $C^{l,\beta}_\oploc$ tensors over such manifolds for $l+\beta\leqslant (k-1)+\alpha$.
\newhead{Cheeger/Gromov Convergence of Abstract Manifolds}
\newsubhead{Cheeger/Gromov Convergence}
\noindent In this appendix, we prove the classical compactness theorem of Riemannian geometry in a form which is most
appropriate for our uses.
\medskip
\noindent We begin by making the following definition:
\proclaim{Definition \nextprocno}
\noindent Let $(M,p)$ be a pointed Riemannian manifold, and let $g$ be the Riemannian metric over $M$.
\medskip
\noindent Let $k\in\Bbb{N}$ be a positive integer, and let $\alpha$ be a real number in $(0,1]$. Let $K,\rho,R>0$ be positive real numbers.
A {\emph $(K,\rho)$-optimal $C^{k,\alpha}$ atlas of $(M,p)$ over a radius $R$} is a family $\Cal{A}=(x_q,\Omega_q,B_\rho(0))_{q\in{B_R(p)}}$
of $C^{k,\alpha}$ coordinate charts of $M$ such that:
\medskip
\myitem{(1)} for all $q\in B_R(p)$, if $A$ is the matrix representation of $(x_q)_*g$ with respect to the Euclidean metric over
$B_\rho(0)$, then:\par
$$\matrix
\|A\|_{C^{k-1,\alpha}(B_\rho(0))} \hfill\leqslant K,\hfill\cr
\|A^{-1}\|_{C^{k-1,\alpha}(B_\rho(0))} \hfill\leqslant K,\hfill\cr
\endmatrix$$
\noindent and,
\medskip
\myitem{(2)} for all $q,q'\in B_R(p)$, and for every ball $B$ contained in $x_q(\Omega_q\minter\Omega_{q'})$:\par
$$
\|x_{q'}\circ x_q^{-1}\|_{C^{k,\alpha}(B)} \leqslant K.
$$
\noindent For all $q\in B_R(p)$, we will refer to $(x_q,\Omega_q,B_\rho(0))$ abusively as a {\emph $(K,\rho)$-optimal $C^{k,\alpha}$
chart of $M$ about $q$}.
\medskip
\noindent We say that $(M,p)$ is {\emph $(K,\rho)$-optimisable over a radius $R$} if such an atlas exists.
\endproclaim
\proclabel{DefnOptimalAtlas}
\noindent Let $\Cal{A}$ be such an atlas of $M$. Let $d$ be the metric (distance structure) generated over $M$ by $g$. For all
$q\in B_R(p)$, let $D_q$ be the distance function generated over $B_\rho(0)$ by $(x_q)_*g$. 
For every $x\in B_\rho(0)$, there exists $\rho'<\rho$ which only
depends on $\rho$, $\left|x\right|$ and $K$ such that, if we denote $\Omega'_q=x^{-1}_q(B_{\rho'}(x))$, then $x_q$ defines an isometry between
$(\Omega'_q,d)$ and $(B_{\rho'}(x),D_q)$. For example, if we choose $x=0$, and if we choose 
$\rho'$ such that:
$$
\rho'\leqslant\frac{\rho}{2K+1},
$$
\noindent then, we find that, for all $y\in B_{\rho'}(0)$:
$$
2D_q(y,0) < D_q(y,\partial B_\rho(0)).
$$
\noindent In this case, if $y,y'\in B_\rho(0)$, then the shortest curve in $M$ joining $x_q^{-1}(y)$ to $x_q^{-1}(y')$ remains within
$\Omega_q$, and thus:
$$
D_q(y,y') = d(x^{-1}_q(y),x^{-1}_q(y')).
$$
\proclaim{Definition \nextprocno}
\noindent We refer to $\rho'$ as the {\emph isometric radius of $\rho$ about $\left|x\right|$ with resect to $K$}.
\endproclaim
\noindent Let $(M_n,p_n)_\ninn$ be a sequence of complete Riemannian manifolds. For all $n$, let $g_n$ be the Riemannian metric
on $M_n$ and let $d_n$ be the metric (distance structure) generated over $M_n$ by $g_n$.
\medskip
\noindent Let $k\in\Bbb{N}$ be a positive whole number and let $\alpha$ be a real number in $(0,1)$. We assume that there exists:
\medskip
\myitem{(1)} a sequence $\suite{R}{n}$ of positive real numbers such that $\suite{R}{n}\uparrow\infty$,
\medskip
\myitem{(2)} a sequence $\suite{N}{n}$ of positive integers such that $\suite{N}{n}\uparrow\infty$,
\medskip
\myitem{(3)} a sequence $\suite{\epsilon}{n}$ of positive real numbers, and
\medskip
\myitem{(4)} a sequence $\suite{K}{n}$ of positive real numbers,
\medskip
\noindent such that, for all $n$, for all $m\geqslant N_n$, there exists a $(K_n,\epsilon_n)$-optimal $C^{k,\alpha}$ atlas $\Cal{A}_{m,n}$
of $(M_m,p_m)$ over a radius $R_n+1$.
\medskip
\noindent We obtain the following result:
\proclaim{Theorem \nextprocno}
\noindent There exists a complete pointed Riemannian manifold $(M_0,p_0)$ of type $C^{k,\alpha}$ such that $(M_n,p_n)_\ninn$ converges to $(M_0,p_0)$ in the pointed weak $C^{k,\alpha}$ Cheeger/Gromov topology.
\medskip
\noindent Moreover, for all $n$, there exists a $(K_n,\epsilon_n)$-optimal $C^{k,\alpha}$ atlas 
$\Cal{A}_{0,n}$ of 
$(M_0,p_0)$ over a radius $R_n+1$.
\medskip
\noindent Finally, there exists a sequence $(\varphi_n)_\ninn$ of weak $C^{k,\alpha}$ convergence mappings of 
$(M_n,p_n)_\ninn$ with respect to $(M_0,p_0)$ such that, for all $m\in\Bbb{N}$, if:
\medskip
\myitem{(1)} $(q_n)_\ninn\in(M_n)_\ninn$ and $q_0\in M$ are such that $(q_n)_\ninn$ converges to $q_0$,
\medskip
\myitem{(2)} for all $n\in\Bbb{N}\munion\left\{0\right\}$, $(x_{q_n},\Omega_{q_n},B_{\epsilon_n}(0))$ is a
$(K_m,\epsilon_m)$-optimal $C^{k,\alpha}$ chart about $q_n$, and
\medskip
\myitem{(3)} $\epsilon_m'$ is the isometric radius of $\epsilon_m$ with respect to $K_m$,
\medskip
\noindent then, there exists a (distance preserving) $C^{k,\alpha}$ mapping
$\alpha:B_{\epsilon_m'}(0)\rightarrow B_{\epsilon_m'}(0)$ such that, after extraction of a subsequence, 
$(x_{q_n}\circ\varphi_n\circ x_{q_0}^{-1})_\ninn$ converges to $\alpha$ in the weak $C^{k,\alpha}$-topology. In otherwords, 
for every compact subset $K$ of  $B_{\epsilon'_m}(0)$, there exists $M\in\Bbb{N}$ such that:
\medskip
\myitem{(1)} for all $n\geqslant M$, $(x_{q_n}\circ\varphi_n\circ x_{q_0}^{-1})$ is defined over $K$ and its restriction
to $K$ is a diffeomorphism onto its image, and
\medskip
\myitem{(2)} for all $\beta<\alpha$:
$$
(\|x_{q_n}\circ\varphi_n\circ x_{q_0}^{-1}-\alpha\|_{C^{k,\beta}(K)})_{n\geqslant M}\rightarrow 0.
$$
\endproclaim
\proclabel{ThmConvergenceOfRiemannianGeometry}
\remark We will refer to such a sequence of convergence mappings as a sequence of {\emph optimal convergence mappings
of $(M_0,p_0)$ with respect to $(M_n,p_n)_\ninn$}.
\medskip
\noindent We obtain this result in many steps. We begin by constructing the limit as a metric space:
\proclaim{Proposition \nextprocno}
\noindent There exists a sequence $(\delta_n)_\ninn\in\Bbb{R}^+$ and a sequence of functions
$(\opCov_n)_\ninn:(0,\delta_n)\rightarrow\Bbb{N}$ such that, for all $n$, for all $\delta<\delta_n$, and for all
$m\geqslant N_n$, the ball $\overline{B_{R_n+1}(p_m)}$ may be covered by $\opCov_n(\delta)$ balls of radius
$\delta$.
\endproclaim
\remark This result is proved by bounding the volumes of small balls from above and below.
\medskip
\proof A proof of this proposition may be found in \cite{Peterson}.\qed
\medskip
\noindent This yields:
\proclaim{Proposition \nextprocno}
\noindent There exists a family $(M_{0,n},d_{0,n},p_{0,n})_\ninn$ of pointed compact metric spaces such that, for all
$n$:
$$
(\overline{B_{R_n+1}(p_m)})_{m\geqslant N_n}\rightarrow (M_{0,n},d_{0,n},p_{0,n}),
$$
\noindent in the Gromov/Haussdorf topology.
\endproclaim
\proof A proof of this proposition may be found in \cite{Peterson}.\qed
\medskip
\noindent By uniqueness of convergence in the Gromov/Haussdorf topology, for all $n\geqslant n'$, there exists a
distance preserving map, sending $(M_{0,n'},p_{0,n'})$ into the ball of radius $R_{n'+1}$ about $p_{0,n'}$ in
$(M_{0,n},p_{0,n})$. Moreover, this map is unique up to isometries of the domain. Consequently, we may take the union of these limiting spaces to obtain:
\proclaim{Proposition \nextprocno}
\noindent There exists a pointed locally compact metric space $(M_0,d_0,p_0)$ such that $(M_n,d_n,p_n)_\ninn$ converges to
$(M_0,d_0,p_0)$ in the pointed Gromov/Haussdorf topology.
\endproclaim
\remark This space is seperable since it is a union of countably many compact sets.
\medskip
\remark We define the set $X$ by:
$$
X = M_0 \munion (\munion_{n\in\Bbb{N}} M_n).
$$
\noindent By definition of the Gromov/Haussdorf topology, we may suppose that there exists a complete metric 
(distance structure) $d$ over $X$ which coincides over $M_n$ with $d_n$ for every $n\in\Bbb{N}\munion\left\{0\right\}$.
We thus adopt the convention that if $(q_n)_\ninn\in(M_n)_\ninn$ is a sequence and if $q_0\in M_0$, then the sequence
$(q_n)_\ninn$ converges to $q_0$ if and only if it converges to $q_0$ in $X$ with respect to the metric $d$.
\medskip
\noindent Since $M_0$ is seperable, it contains a countable dense subset which we will denote by $Q$. For every $q\in Q$,
we define $k_q$ by:
$$
k_q = \minf\left\{k\in\Bbb{N}\text{ s.t. } q\in B_{R_k+1}(p_0)\right\}.
$$
\noindent For all $q$, let $\epsilon_{k_q}'$ be the isometric radius of $\epsilon_{k_q}$ with respect to $K_{k_q}$. Since $Q$ is dense, we have:
$$
M_0 = \munion_{q\in Q} B_{\epsilon_{k_q}'}(q).
$$
\noindent For every $q\in Q$, we define $(q_n)_\ninn\in (M_n)_\ninn$ to be a sequence which converges towards $q$. For all
sufficiently large $n$, we have:
$$
q_n \in B_{R_{k_q}+1}(p_n).
$$
\noindent These definitions will be of use to us in the sequel.
\medskip
\noindent We now furnish the metric space underlying the limiting manifold with a $C^{k,\alpha}$ differential structure 
and a $C^{k-1,\alpha}$ Riemannian metric:
\proclaim{Proposition \nextprocno}
\noindent There exists a canonical maximal $C^{k,\alpha}$ atlas $\Cal{A}$ and a canonical $C^{k-1,\alpha}$ Riemannian metric $g_0$
over $M_0$ such that the metric (distance structure) generated by $g_0$ over $M_0$ coincides with $d_0$.
\medskip
\noindent Moreover, for all $n$, the pointed manifold $(M_0,p_0)$ is $(K_n,\epsilon_n)$-normalisable over a radius $R_n$.

\endproclaim
\proclabel{PropDiffPlusRiemStruct}
\proof Without loss of generality, we may suppose that, for all $n$:
$$
q_n\in B_{R_{k_q}+1}(p_n).
$$
\noindent For all $n\geqslant N_{k_q}$, let $(x_{q,n},\Omega_{q,n},B_{\epsilon_{k_q}}(0))$ be a 
$(K_q,\epsilon_{k_q})$-optimal $C^{k,\alpha}$ chart of $M_n$ about $q_n$.
\medskip
\noindent For all $n$, we define the metric $g_{q,n}$ over $B_{\epsilon_{k_q}}(0)$ by:
$$
g_{q,n} = (x_{q,n})_*g_n.
$$
\noindent By the classical Arzela-Ascoli theorem, after extraction of a subsequence, there exists a $C^{k-1,\alpha}$ metric
$g_q$ over $B_{\epsilon_{k_q}}(0)$ such that $(g_{q,n})_\ninn$ converges to $g_q$ in the weak $C^{k-1,\alpha}$ topology. Since the metrics
$(g_{q,n})_\ninn$ are uniformly bounded below, it follows that $g_q$ is positive definite.
\medskip
\noindent For all $n$, let $d_{q,n}$ be the
distance structure generated over $B_{\epsilon_{k_q}}(0)$ by $g_{q,n}$. Let $d_q$ be the distance structure generated by $g_q$.
\medskip
\noindent Since, for all $n$, $(x_{q,n})^{-1}$ is locally distance preserving, by the classical
Arzela-Ascoli theorem, there exists a locally distance preserving mapping $\xi_q:B_{\epsilon'_{k_q}}(0)\rightarrow M_0$ such that
$(x^{-1}_{q,n})_\ninn$ converges locally uniformly to $\xi_q$.
\medskip
\noindent Let $x$ be a point in $B_\epsilon(0)$. Let $\rho$ be the isometric radius of 
$\epsilon_{k,q}$ about $\left|x\right|$ with respect to $K$. Since the restriction of
$\xi_q$ to $B_\rho(x)$ preserves distance, it is a homeomorphism onto its image. Moreover, the image of $\xi_q$ is an
open set. Indeed, let $y$ be a point in $B_{\rho}(x)$. Let $\delta$ be such that the closed ball of radius $\delta$ about $y$
with respect to $d_q$ is a compact subset of $B_{\rho}(x)$. We find that the image of the interior of this ball under the action of
$\xi_q$ is precisely the open ball of radius $\delta$ about $\xi_q(y)$ in $M_0$, and the openness of $\xi_q$ now follows. Since, for every $n$, $(x_{q,n})^{-1}$ is a homeomorphism, it follows that
$\xi_q$, being the locally uniform limit of a sequence of homeomorphisms, is also a homeomorphism
(see, for example, lemma $2.2.2$ of \cite{SmiE}). We thus define $x_q=\xi_q^{-1}$ and 
$\Omega_q=\xi_q(B_{\epsilon_{k_q}'}(0))$ and we obtain a chart $(x_q,\Omega_q,B_{\epsilon'_{k_q}}(0))$ of $M_0$ about $q$.
\medskip
\noindent Let $q'$ be another point in $Q$. We construct $(x_{q'},\Omega_{q'},B_{\epsilon_{k_{q'}}}(0))$ as for $q$. We suppose that
$\Omega_q\minter\Omega_{q'}\neq\emptyset$. Let $K$ be a compact subset of $\Omega_q\minter\Omega_{q'}$. We have
$x_q(K)\subseteq B_{\epsilon_{k_q}}(0)$. For all sufficiently large $n$, we have:
$$
x^{-1}_{q_n}(x_q(K)) \subseteq \Omega_{q_n'}.
$$
\noindent Indeed, let $y$ be a point in $x_q(K)$. Let us define $z\in B_{\epsilon_{q'}}(0)$ by:
$$
z = (x_{q'}\circ x_q^{-1})(y).
$$
\noindent Let $\epsilon_1$ be such that $B_{\epsilon_1}(z)\subseteq B_{\epsilon'_{k_{q'}}}(0)$. Since the mappings
$(x_{q_n'})_{n\geqslant N_{q'}}$ are uniformly bilipschitzian, there exists $\epsilon_2$ such that, for all
$n\geqslant N_{q'}$:
$$\matrix
B_{\epsilon_2}(x^{-1}_{q_n'}(z)) \hfill&\subseteq x^{-1}_{q_n'}(B_{\epsilon_1}(z))\hfill\cr
&\subseteq \Omega_{q_n'}.\hfill\cr
\endmatrix$$
\noindent Since $(x_{q_n}^{-1}(y))_\ninn$ and $(x^{-1}_{q_n'}(z))_\ninn$ both converge towards $x_q^{-1}(y)=x_{q'}^{-1}(z)$, it
follows that, for all sufficiently large $n$:
$$
(x^{-1}_{q_n}(y)) \subseteq B_{\epsilon_2/2}(x_{q_n'}^{-1}(z)).
$$
\noindent Consequently, there exists $\epsilon_3$ such that, for all sufficiently large $n$:
$$\matrix
x_{q_n}^{-1}(B_{\epsilon_3}(y)) \hfill&\subseteq B_{\epsilon_2/2}(x_{q_n}^{-1}(y))\hfill\cr
&\subseteq B_{\epsilon_2}(x_{q_n'}^{-1}(z))\hfill\cr
&\subseteq \Omega_{q_n'}.\hfill\cr
\endmatrix$$
\noindent Since $y$ is arbitrary, and since $K$ is compact, it follows that for sufficiently large $n$:
$$
x^{-1}_{q_n}(x_q(K)) \subseteq \Omega_{q_n'},
$$
\noindent and we obtain the desired result.
\medskip
\noindent Let $B$ be an open ball whose closure is contained in 
$x_q(\Omega_q\minter\Omega_{q'})\subseteq B_{\epsilon_q}(0)$. In particular, $x_q^{-1}(\overline{B})$
is compact, and thus, for sufficiently large $n$, $x_{q_n'}\circ x_{q_n}^{-1}$ is defined over $B$. By our 
hypotheses, for sufficiently large $n$:
$$
\|x_{q_n'}\circ x_{q_n}^{-1}\|_{C^{k,\alpha}(B)}\leqslant K.
$$
\noindent Consequently, by the classical Arzela-Ascoli theorem, there exists $\varphi\in C^{k,\alpha}(B)$ such that, after
extraction of a subsequence, $(x_{q_n'}\circ x_{q_n}^{-1})_\ninn$ converges to $\varphi$ in the weak 
$C^{k,\alpha}$ topology over $B$. However, for all sufficiently large $n$:
$$
x_{q_n'}^{-1}\circ(x_{q_n'}\circ x_{q_n}^{-1}) = x_{q_n}^{-1}.
$$
\noindent Consequently, by taking limits, we obtain:
$$\matrix
&\xi_{q'}\circ\alpha \hfill&=\xi_q\hfill\cr
\Rightarrow\hfill&\alpha\hfill&= x_{q'}\circ x_q^{-1}.\hfill\cr
\endmatrix$$
\noindent In particular, it follows that $(x_{q'}\circ x_q^{-1})$ is of type $C^{k,\alpha}_\oploc$.
\medskip
\noindent It thus follows that $(x_q,\Omega_q,B_{\epsilon_{k_q}}(0))_{q\in Q}$ is a $C^{k,\alpha}$ atlas over $M_0$.
\medskip
\noindent For all $q,q'\in Q$, and for all sufficiently large $n$:
$$
((x_{q_n}')_*g_n)|_{x_{q_n'}(\Omega_{q_n'}\minter\Omega_{q_n})} =
(x_{q_n'}\circ x_{q_n}^{-1})_*((x_{q_n})_*g_n)|_{x_{q_n}(\Omega_{q_n'}\minter\Omega_{q_n})}.
$$
\noindent Thus, taking limits, we obtain:
$$
g_{q'}|_{x_q'(\Omega_{q'}\minter\Omega_q)} = (x_{q'}\circ x_q^{-1})_* g_q|_{(x_q(\Omega_{q'}\minter\Omega_q)}.
$$
\noindent It follows that the family $(g_q)_{q\in Q}$ defines a $C^{k-1,\alpha}$ Riemannian metric $g_0$ over $M_0$.
Moreover, with these differential and Riemannian structures, for all $n\in\Bbb{N}$, $(M_0,p_0)$ is 
$(K_n,\epsilon_n)$-normalisable over a radius of $R_n$.
\medskip
\noindent Let $p_0$,$p_1$ be two points in $M_0$. Let $\epsilon\in\Bbb{R}^+$ be a positive real number. There exist
sequences $(p_{0,n})_\ninn,(p_{1,n})_\ninn\in (M_n)_\ninn$ such that $(p_{0,n})_\ninn$ and $(p_{1,n})_\ninn$ converge
respectively to $p_0$ and to $p_1$. For sufficiently large $n$, we may suppose that:
$$
d_n(p_{0,n},p_{1,n}) \leqslant d_0(p_0,p_1) + \epsilon/2.
$$
\noindent For all $n$, there exists a continuous path $\gamma_n:I\rightarrow M_n$ joining $p_{0,n}$ to $p_{1,n}$ such that:
$$
\opLength(\gamma_n) \leqslant d_n(p_{0,n},p_{q,n})+\epsilon/2.
$$
\noindent Thus, for sufficiently large $n$:
$$
\opLength(\gamma_n) \leqslant d_n(p_0,p_1) + \epsilon.
$$
\noindent We may suppose that every $\gamma_n$ is parametrised by a constant factor of arc length. The classical
Arzela-Ascoli theorem now tells us that there exists $\gamma_0:I\rightarrow M_0$ such that $\suite{\gamma}{n}$
converges uniformly to $\gamma_0$ over $I$. In particular:
$$
\gamma_0(0) = p_0,\qquad \gamma_0(1) = p_1,
$$
\noindent and:
$$
\opLength(\gamma_0) \leqslant d_0(p_0,p_1) + \epsilon.
$$
\noindent Since $\epsilon\in\Bbb{R}^+$ is arbitrary, it follows that $d_0$ is a length metric.
\medskip
\noindent However, for all $q$, if $\epsilon_{k_q}'$ is the isometric radius of $\epsilon_{k_q}$ about $0$ with respect to $K$, then $x_q:(\Omega_q,d_0)\rightarrow(B_{\epsilon'_{k_q}}(0),d_q)$ is an isometry of metric 
spaces. Thus $d_0$ coincides everywhere locally with the distance structure generated over $M_0$ by $g_0$. Since $d_0$ 
is a length metric, the result now follows.\qed
\medskip
\goodbreak
\noindent For all $q\in Q$, we construct $(x_q,\Omega_q,B_{\epsilon_{k_q}(0)})$ as in the proof of proposition 
\procref{PropDiffPlusRiemStruct}. We now make the following definition:
\proclaim{Definition \nextprocno}
\noindent Let $U$ be an open subset of $M_0$ containing $p_0$. Let $k\in\Bbb{N}$ be a positive integer and let $\alpha$
be a real number in $(0,1]$. A sequence of {\emph $C^{k,\alpha}$ strong convergence mappings of $(M_n,p_n)_\ninn$ over $U$ 
with respect to $Q$} is a sequence $(\varphi_n)_\ninn$ such that:
\medskip
\myitem{(1)} for all $n$, $\varphi_n:(U,p_0)\rightarrow (M_n,p_0)$ is a $C^{k-1,\alpha}_\oploc$ diffeomorphism onto its image, and
\medskip
\myitem{(2)} for all $q,q'\in Q$, the sequence $(x_{q_n'}\circ\varphi_n\circ x_q^{-1})_\ninn$ converges to
$(x_{q'}\circ x_q^{-1})$ in the weak $C^{k,\alpha}_\oploc$ topology. In other words, for every compact set
$K\subseteq x_q(\Omega_q\minter\Omega_{q'}\minter U)\subseteq B_{\epsilon_{k_q}'}(0)$:
\medskip
\myitem{(a)} there exists $N$ such that, for all $n\geqslant N$:
$$\matrix
x_q(K) \hfill&\subseteq \Omega_q,\hfill\cr
\varphi_n(x_q(K)) \hfill&\subseteq \Omega_{q_n'}.\hfill\cr
\endmatrix$$
\noindent and,
\medskip
\myitem{(b)} $(x_{q_n'}\circ\varphi_n\circ x_q^{-1})_{n\geqslant N}$ converges to $x_{q'}\circ x_q^{-1}$ in the weak
$C^{k,\alpha}$ topology over $K$.
\endproclaim
\remark By taking subsequences, we may suppose that, for all $q\in Q$, $(x_{q_n}^{-1}\circ x_q)$ is a sequence of $C^{k,\alpha}$ strong convergence mappings of $(M_n,p_n)_\ninn$ over $U$ with
respect to $Q$.
\medskip
\noindent We now prove theorem \procref{ThmConvergenceOfRiemannianGeometry} by induction. The induction step is
guaranteed by the following result:
\proclaim{Proposition \nextprocno}
\noindent Let $U$ be an open subset of $M$ containing $p_0$. Let $(\varphi_n)_\ninn$ be a sequence of strong convergence 
mappings of $(M_n,p_n)$ over $U$ with respect to $Q_0$. For every $q\in Q$ and for
every relatively compact open subset $V\subseteq U\munion \Omega_q$, after extraction of a subsequence of 
$(M_n,p_n)_\ninn$, there exists a sequence $(\psi_n)_\ninn$ of strong convergence mappings of $(M_n,p_n)_\ninn$ over $V$ with respect to $Q$.
\endproclaim
\proclabel{PropConvergenceMappings}
\proof Let $(\phi_U,\phi_q)$ be a $C^{k,\alpha}$-partition of unity of $U\munion\Omega_q$ subordinate to the cover
$(U,\Omega_q)$. 
\medskip
\noindent Define $K\subseteq U\minter\Omega_q$ by:
$$
K = \overline{V}\minter\opSupp(\phi_q)\minter\opSupp(\phi_U).
$$
\noindent Since $(\varphi_n)_{n\geqslant N}$ is a sequence of strong convergence mappings, there exists $M_1$ such that
for $n\geqslant M_1$:
$$
\varphi_n(K) \subseteq \Omega_{q_n}.
$$
\noindent For $n\geqslant M_1$, we define $\psi_n:V\rightarrow M_n$ by:
$$
\psi_n(p) = \left\{
\matrix
\varphi_n(p) \hfill&\text{ if }p\in\opSupp(\phi_q)^C\minter V,\hfill\cr
(x_{q_n}^{-1}\circ x_q)(p)\hfill&\text{ if }p\in\opSupp(\phi_U)^C\minter V,\hfill\cr
x_{q_n}^{-1}(\phi_U(p)(x_{q_n}\circ\varphi_n)(p) + \phi_q(p)x_q(p))\hfill&\text{ otherwise (i.e., if $p\in K$).}\hfill\cr
\endmatrix\right.
$$
\noindent For all $n$, since $\varphi_n$ is of type $C^{k,\alpha}_\oploc$, so is $\psi_n$.
\medskip
\noindent Let $p$ be a point in $\opSupp(\phi_q)^C\munion\opSupp(\phi_U)^C$. Let $\Omega$ be a neighbourhood of $p$ in\break
$\opSupp(\phi_q)^C\munion\opSupp(\phi_U)^C$. Since, for all $n$, $\psi_n$ coincides over $\Omega$ either with
$\varphi_n$ or with $(x_{q_n}^{-1}\circ x_q)$, it follows trivially that $(\psi_n)_\ninn$ is a sequence of
strong convergence mappings of $(M_n,p_n)_\ninn$ over $\Omega$ with respect to $Q$. It thus remains to study what
happens near points in $K$.
\medskip
\noindent Let $p$ be a point in $K$. Let $\Omega$ be a relatively compact neighbourhood of $p$ in $\Omega_q$. There
exists $N\geqslant N_{k_q}$ such that, for all $n\geqslant N$:
$$
\Omega\subseteq \Omega_q\minter\varphi_n^{-1}(\Omega_{q_n}).
$$
\noindent For $n\geqslant N$, and for $z\in x_q(\Omega)$, we have:
$$
(x_{q_n}\circ\psi_n\circ x_q^{-1})(z) = (\phi_U\circ x_q^{-1})(z)(x_{q_n}\circ\varphi_n\circ x_q^{-1})(z) + (\phi_q\circ x_q^{-1})(z)z.
$$
\noindent Since $(\varphi_n)_\ninn$ is a sequence of $C^{k,\alpha}$ strong convergence mappings over $U$ with respect to $Q$, it follows that, since $(x_{q_n}\circ\varphi_n\circ x_q^{-1})(z)_\ninn$ converges to $x_q\circ x_q^{-1}=\opId$ in the weak
$C^{k,\alpha}_\oploc$ topology over $x_q(\Omega)$. Consequently $(x_{q_n}\circ\psi_n\circ x_q^
{-1})(z)_\ninn$ also converges to the identity in the weak $C^{k,\alpha}_\oploc$ topology
over $x_q(\Omega)$.
\medskip
\noindent Let $q'$,$q''$ be two other points in $Q$. Let us suppose that $p\in\Omega_q\minter\Omega_{q'}\minter\Omega_{q''}$
and let $\hat{\Omega}$ be a relatively compact neighbourhood of $p$ in $\Omega_q\minter\Omega_{q'}\minter\Omega_{q''}$.
\medskip
\noindent For sufficiently large $n$, $x_{q_n''}\circ\psi_n\circ x_{q'}^{-1}$ is defined over
$\hat{\Omega}$, and, for all $z\in\hat{\Omega}$:
$$
(x_{q_n''}\circ\psi_n\circ x_{q'}^{-1})(z) = (x_{q_n''}\circ x_{q_n}^{-1})\circ (x_{q_n}\circ\psi_n\circ x_q^{-1})\circ(x_q\circ x_{q'}^{-1}).
$$
\noindent It thus follows that $(x_{q_n''}\circ\psi_n\circ x_{q'}^{-1})_\ninn$ converges to 
$(x_{q''}\circ x_{q'}^{-1})$ over $\hat{\Omega}$ in the weak $C^{k,\alpha}_\oploc$ topology.
\medskip
\noindent It thus follows that for all $q',q''\in Q$, $(x_{q_n''}\circ\psi_n\circ x_{q'}^{-1})_\ninn$ converges to
$(x_{q''}\circ x_{q'}^{-1})$ in the weak $C^{k,\alpha}$ topology over $x_{q'}(\Omega_{q'}\minter\Omega_{q''}\minter V)$. It thus
follows that, if we can show that the, for sufficiently large $n$, the restriction of $\psi_n$ to $V$ is a 
diffeomorphism onto its image, then $(\psi_n)_\ninn$ will be a sequence of strong convergence mappings of $(M_n,p_n)$ over 
$V$ with respect to $Q$.
\medskip
\noindent For all $q',q''\in Q$, for sufficiently large $n$, $x_{q_n}''\circ\psi_n\circ x_{q'}^{-1}$ is
everywhere a local diffeomorphism. Consequently, for all sufficiently large $n$, $\psi_n$ is everywhere a local
diffeomorphism. It now remains to show that, for all sufficiently large $n$, $\psi_n$ is injective.
\medskip
\noindent Suppose that there exists $(x_n)_\ninn,(y_n)_\ninn\in V$ such that, for all
$n\in\Bbb{N}$:
$$
x_n\neq y_n,\qquad \psi_n(x_n) = \psi_n(y_n).
$$
\noindent By compactness, we may assume that there exists $x_0,y_0\in\overline{V}$ such that $(x_n)_\ninn$ and
$(y_n)_\ninn$ converge respectively to $x_0$ and $y_0$. Since, for all $q',q''\in Q$, the sequence
$(x_{q_n''}\circ\psi_n\circ x_{q'}^{-1})_\ninn$ converges to $x_{q''}\circ x_{q'}^{-1}$, we find that $x_0=y_0$. 
Since this sequence converges in the weak $C^{k,\alpha}_\oploc$ topology and since $k\geqslant 1$, it follows that 
$x_n=y_n$ for all sufficiently large $n$ (see, for example, lemma $2.2.3$ of \cite{SmiE}). This is absurd, and it follows that there exists $N$ such that for 
$n\geqslant N$, the application $\psi_n:V\rightarrow M_n$ is a diffeomorphism
onto its image. The result now follows.\qed
\medskip
\noindent We may now prove theorem \procref{ThmConvergenceOfRiemannianGeometry}:
\medskip
{\sl\noindent Proof of theorem \procref{ThmConvergenceOfRiemannianGeometry}:\ } By proposition \procref{PropConvergenceMappings}, for every
finite family $Q'\subseteq Q$, and for every relatively compact subset $V$ of $\munion_{q\in Q'}\Omega_q$, we may find a positive integer $N$
and a sequence of $C^{k,\alpha}$ strong convergence mappings $(\varphi_n)_{n\geqslant N}$ over $V$. Using a diagonal argument (\`a la Cantor), we may
thus find:
\medskip
\myitem{(1)} a sequence of real numbers $(\rho_n)_\ninn$ such that $(\rho_n)_\ninn\uparrow\infty$, and
\medskip
\myitem{(2)} for all $n\in\Bbb{N}$, a local diffeomorphism $\varphi_n:(B_{\rho_n}(p_0),p_0)\rightarrow (M_n,p_0)$,
\medskip
\noindent such that, for all $N\in\Bbb{N}$, the sequence $(\varphi_n)_{n\geqslant N}$ is a sequence of strong
convergence mappings over $B_{\rho_N}(p_0)$.
\medskip
\noindent Consequently, for all $N\in\Bbb{N}$, and for all $n\geqslant N$, the restriction of $\varphi_n$ to
$B_{\rho_n}(p_0)$ is a diffeomorphism onto its image. Moreover, since $(\varphi_n)_{n\geqslant N}$ is a sequence of
strong convergence mappings over $B_{\rho_N}(p_0)$, we find that $(\varphi_n^*g_n)_{n\geqslant N}$ converges to $g_0$
in the weak $C^{k,\alpha}_\oploc$ topology over $B_{\rho_N}(p_0)$ and the first result now follows.
\medskip
\noindent The second result follows directly from proposition \procref{PropDiffPlusRiemStruct} and the third result is an 
immediate consequence of the fact that $(\varphi_n)_\ninn$ is a sequence of strong convergence mappings.\qed
\newhead{Bibliography}
{\leftskip = 5ex \parindent = -5ex
\leavevmode\hbox to 4ex{\hfil\cite{Cheeger}}\hskip 1ex{Cheeger J., Finiteness theorems for Riemannian manifolds, {\sl Amer. J. Math.} {\bf 92} (1970), 61--74}
\medskip
\leavevmode\hbox to 4ex{\hfil\cite{Corlette}}\hskip 1ex{Corlette K., Immersions with bounded curvature, {\sl Geom. Dedicata} {\bf 33} (1990), no. 2, 153--161}
\medskip
\leavevmode\hbox to 4ex{\hfil\cite{GromA}}\hskip 1ex{Gromov M., {\sl Metric Structures for Riemannian and Non-Riemannian Spaces}, Progress in Mathematics, {\bf 152}, Birkh\"auser, Boston, (1998)}
\medskip
\leavevmode\hbox to 4ex{\hfil\cite{Peterson}}\hskip 1ex{Peterson P., {\sl Riemannian Geometry}, Graduate Texts in Mathematics, {\bf 171}, Springer Verlag, New York, (1998)}
\medskip
\leavevmode\hbox to 4ex{\hfil\cite{SmiA}}\hskip 1ex{Smith G., Special Legendrian structures and Weingarten problems, Preprint, Orsay (2005)}%
\medskip
\leavevmode\hbox to 4ex{\hfil\cite{SmiE}}\hskip 1ex{Smith G., Th\`ese de doctorat, Paris (2004)}%
\par
}%
\enddocument